\documentclass[11pt, twoside]{article}

\usepackage{amsmath}
\usepackage{amsfonts}
\usepackage{amssymb}
\usepackage{enumerate}

\usepackage{array}
\usepackage{epic}
\usepackage{eepic}
\usepackage{color}

\usepackage{graphicx}
\usepackage{fancyhdr}
\usepackage{geometry}
\begin{document}

\newtheorem{theorem}{Theorem}[section]
\newtheorem{lemma}[theorem]{Lemma}
\newtheorem{corollary}[theorem]{Corollary}
\newtheorem{proposition}[theorem]{Proposition}
\newtheorem{example}[theorem]{Example}
\newtheorem{examples}[theorem]{Examples}
\newtheorem{definition}[theorem]{Definition}
\newtheorem{remarks}[theorem]{Remarks}
\newtheorem{exercise}[theorem]{Exercises}
\def\endbox{\begin{flushright}$\Box$\end{flushright}}
\def\bproof{{\noindent\bf Proof: }}
\def\eproof{\begin{flushright}$\dashv$\end{flushright}}
\def\<{\langle}
\def\>{\rangle}
\def \into#1#2#3{#1 \buildrel #2 \over \hookrightarrow #3}
\def \con#1{{\text {\rm Con}} (#1)}
\def \rpow#1#2#3{\left. {#1} ^ {#2} \right/ {#3}}

\pagestyle{fancy} \fancyhf{} \title{Enlargement of Filtrations --- A
Primer}
\author{P. Ouwehand\\ \\African Institute for Financial Markets and Risk Management\\University of Cape Town}
\date{}
\maketitle
\vskip0.3cm

\begin{center}\includegraphics[scale=0.5]{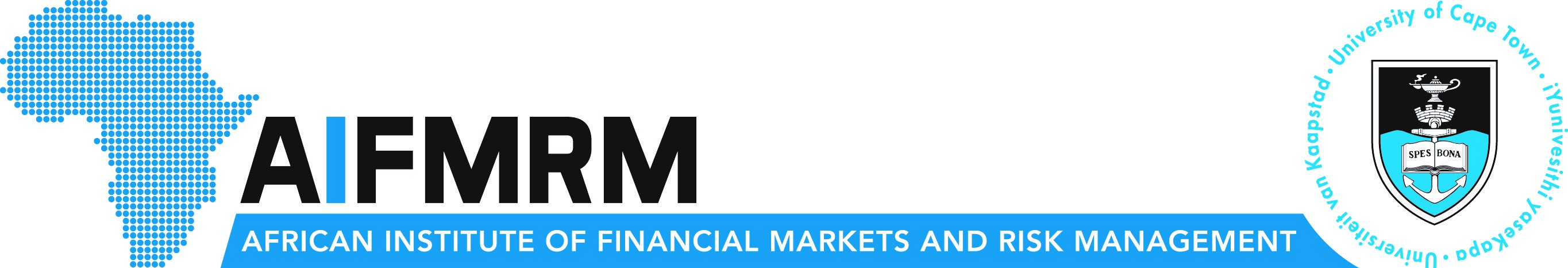}\end{center}

\tableofcontents

 \fancyhead[LE,RO]{\thepage} \fancyhead[LO,RE]{\em
Enlargement of Filtrations}

\newpage

\section{Introduction}
Let $(\Omega,{\cal F},\mathbb P)$ be a probability space equipped
with two filtrations $\mathbb F = ({\cal F}_t)_{t\geq 0}, \mathbb
G = ({\cal G}_t)_{t\geq 0}$, both satisfying the usual conditions,
and such that ${\cal F}_t\subseteq{\cal G}_t$ for all $t\geq 0$.

We will study the following hypothesis:\begin{center}\fbox{
Hypothesis ${\mathbf H'}$: Every $(\mathbb F,\mathbb
P)$--semimartingale is a $(\mathbb G,\mathbb
P)$--semimartingale.}\end{center}
 Note that by the Doob--Meyer
Decomposition Theorem and the Fundamental Theorem of Local
Martingales, $(\mathbf H')$ is equivalent to each of the
following:\begin{enumerate}[(i)]\item Every $(\mathbb F,\mathbb
P)$--local martingale is a $(\mathbb G,\mathbb
P)$--semimartingale. \item Every bounded $(\mathbb F,\mathbb
P)$--martingale is a $(\mathbb G,\mathbb P)$--semimartingale.
\end{enumerate}

\begin{example}\label{example_add_W_T}\rm
Let $W$ be a standard Brownian motion on $(\Omega,{\cal F},\mathbb
P)$, with natural filtration $\mathbb F=({\cal F}_t)_{t\geq 0}$,
and expand $\mathbb F$ by adding knowledge of $W_T$ for some
$T>0$, i.e. consider the enlarged filtration $\mathbb H =({\cal
H}_t)_{t\geq 0}$ which is the smallest filtration satisfying the
usual conditions such that $\sigma(W_T)\cup{\cal F}_t\subseteq
{\cal H}_t$ for all $t\geq 0$. We will show (following
Protter\cite{Protter2004}) that $W$ is a $\mathbb
H$--semimartingale (though no longer a $\mathbb H$--martingale)
with decomposition
\[W_t = M_t+A_t\qquad\text{where}\quad M_t = W_t -\int_0^{t\land
T}\frac{W_T-W_s}{T-s}\;ds\] and $M$ is a $\mathbb H$--martingale.

Without loss of generality, we may assume that $T=1$. To show $M$
is a $\mathbb H$--martingale, we will need to calculate $\mathbb
E[W_t-W_s|{\cal H}_s]$. Since $W_t-W_s$ is independent of ${\cal
F}_s$ for $t>s$, and since ${\cal H}_s$ is (essentially) ${\cal
F}_s\lor \sigma(W_1-W_s)$ it follows that $\mathbb E[W_t-W_s|{\cal
F}_s] = \mathbb E[W_t-W_s|W_1-W_s]$. Symmetry now suggests
that\[\mathbb E[W_t-W_s|W_1-W_s]=\frac{t-s}{1-s}(W_1-W_s)\] To
prove this relation, let $0\leq s<t\leq 1$ be rational numbers.
Then there exists $n,j,k\in\mathbb N$ such that $s=\frac{j}{n},
t=\frac{k}{n}$. For $i=0,1,\dots, n-1$, define $Y_i =
W_{\frac{i+1}{n}}-W_{\frac{i}{n}}$. Then $W_1-W_s
=\sum_{i=j}^{n-1}Y_i, W_t-W_s = \sum_{i=j}^{k-1}Y_i$. Now as the
$Y_i$ are independent and identically distributed, symmetry
dictates that
\[\mathbb E[W_t-W_s|W_1-W_s] =\mathbb
E\left[\sum_{i=j}^{k-1}Y_i\Big|\sum_{i=j}^{n-1}Y_i\right]=
\frac{k-1}{n-j}\sum_{i=j}^{n-1}Y_i = \frac{t-s}{1-s}(W_1-W_s)\]
Now as $W_t-W_s$ is independent of ${\cal F}_s$, we have
\[\mathbb E[W_t-W_s|{\cal H}_s]=\mathbb
E[W_t-W_s|W_1-W_s]=\frac{t-s}{1-s}(W_1-W_s)\] for all rationals
$0\leq s<t\leq 1$. Since $W$ is continuous, a little argument
shows that this relation holds for all $0\leq s<t\leq 1$. Now,
with $M_t$ as defined above, and applying Fubini's Theorem (for
conditional expectations) we have
\[\aligned \mathbb E[M_t-M_s|{\cal H}_s] &= \mathbb E[W_t-W_s|{\cal
H}_s]-\int_s^t\frac{\mathbb E[W_1-W_u|{\cal H}_s]}{1-u}\;du\\
&=\frac{t-s}{1-s}(W_1-W_s)-\int_s^t\frac{(1-u)(W_1-W_s)}{(1-u)(1-s)}\;du\\&=0\endaligned\]
for $0\leq s<t<1$, i.e. $M$ is a $\mathbb H$--martingale on
$[0,1)$. We want to show it that $M$ is a martingale on all of
$\mathbb R^+$. There may be a problem at $t=1$ because of the
possibility of explosions: We require that
$\int_s^1\frac{W_1-W_s}{1-s}\;ds$ remains finite a.s. This will be
guaranteed if we can show that $\mathbb
E[\int_0^1\frac{|W_1-W_s|}{1-s}\;ds]<\infty$. Since the
$L^2$--norm dominates the $L^1$--norm, we see that $\mathbb
E[|W_1-W_s|]\leq \mathbb E[(W_1-W_s)^2]^{\frac12}= \sqrt{1-s}$,
and thus that $\mathbb E[\int_0^1\frac{|W_1-W_s|}{1-s}\;ds]\leq
\int_0^1\frac{1}{\sqrt{1-s}}\;ds<\infty$ as required.

Finally, if $t>1$, then ${\cal F}_t = {\cal H}_t$. It is now easy
to see that $M$ is a $\mathbb H$--martingale.
\endbox
\end{example}

\begin{example}\rm Let $W$ be a standard Brownian motion, with
natural filtration $\mathbb F$, and define $\mathbb G$ by ${\cal
G}_t = {\cal F}_{t+\varepsilon}$, where $\varepsilon >0$. Then
$H'$ fails: $W$ is not a $\mathbb G$--semimartingale.

To prove this, we use the Bichteler--Dellacherie Theorem, i.e. we
show that $W$ is not a good $(\mathbb G,\mathbb P)$--integrator.
It suffices to show that there exists a sequence $H^n$ of $\mathbb
G$--predictable elementary processes such that $H^n\to 0$
uniformly in probability on compacts, yet $(H^n\bullet W)\not\to
0$ ucp (cf. \cite{Protter2004}, Thm. II.11). Let $n_0\in\mathbb N$
be least such that $2^{-n_0}\leq\varepsilon$. For $n\geq N$,
define
\[H^n = \sum_{k=0}^{2^n-1}
\Delta^n_kWI_{(\frac{k}{2^n},\frac{k+1}{2^n}]}\qquad\text{where
}\quad\Delta^n_kW = W_ {\frac{k+1}{2^n}}-W_{\frac{k}{2^n}}\] Note
that if $Y\sim N(0,1)$ and $y\geq 0$, then $\mathbb P(Y\geq y)\leq
\frac{1}{\sqrt{2\pi}y}e^{-\frac12y^2}$, a fact which follows
easily from the identity $e^{-\frac12 y^2}\frac1y\leq
\int_y^\infty e^{-\frac12x^2}(1+\frac{1}{x^2})\;dx$. It follows
that, for $\delta
>0$, we have
\[\mathbb P(|\Delta^n_kW|>\delta)\leq
\frac{2}{2^{n/2}\delta\sqrt{2\pi}}\;e^{-\frac12 2^n\delta^2}\] and
hence that
\[\mathbb P(\sup_{t\leq 1}|H^n_t|>\delta)
\leq\sum_{k=0}^{2^n-1}\mathbb P(|\Delta^n_kW|>\delta) \leq
\frac{2^n\cdot 2}{2^{n/2}\delta\sqrt{2\pi}}\;e^{-\frac12
2^n\delta^2} = \frac{KC}{e^{\frac12 K^2D}}\] where $K=2^{n/2}, C =
\frac{2}{\delta\sqrt{2\pi}}$ and $D = \delta^2$.  It follows
easily that $H^n\to 0$ ucp as $n\to\infty$. However,
\[(H^n\bullet W)_1 = \sum_{k=0}^{2^n-1}(\Delta^n_kW)^2\] so that
$\mathbb E(H^n\bullet W)_1 =1$. Moreover, the family
$\{(H^n\bullet W)_1:n\geq n_0\}$ is clearly ${\cal L}^2$--bounded
(by 5, for example), hence UI. We therefore cannot have
$(H^n\bullet W)\to 0$ ucp.
\endbox
\end{example}
\begin{remarks}\rm\begin{enumerate}[(1.)]\item Protter\cite{Protter2004}, citing Ito,  shows that
Example \ref{example_add_W_T} also applies to L\'evy processes: If
$Z$ is a L\'evy process on $(\Omega,{\cal F},\mathbb P,\mathbb
F)$, and if ${\cal H}$ is the smallest filtration satisfying the
usual conditions such that $\sigma(Z_1)\cup{\cal
F}_t\subseteq{\cal H}_t$ for all $t\geq 0$, then $Z$ is a $\mathbb
H$--semimartingale. In addition, if $\mathbb E|Z_t|<\infty$ for
all $t\geq 0$, then \[M_t = Z_t -\int_0^{t\land
1}\frac{Z_1-Z_s}{1-s}\;ds\] is a $\mathbb H$--martingale. \item
Later, we will generalize Example \ref{example_add_W_T} by
considering enlargements of the form $\sigma(\int_0^\infty
\varphi(s)\;dW_s)$.
\end{enumerate}
\endbox
\end{remarks}

\begin{remarks}\rm
The semimartingale property is closely related to the no--arbitrage
property. Delbaen and Schachermayer\cite{Delbaen_Schachermayer1995}
show that if $S$ is a locally bounded c\`adl\`ag asset price process
satisfying the\lq\lq no free lunch with vanishing risk property"
(NFLVR) for simple portfolios (i.e. linear combinations of
buy--and--hold portfolios $I_{(T_1, T_2]}$, where $T_1\leq T_2$ are
stopping times, then $S$ must be a semimartingale. In other words,
if the semimartingale property is destroyed by enlargement of
filtration, then there must be something very close to arbitrage for
those with access to the enlarged information set.

However, the semimartingale property is not sufficient to guarantee
no--arbitrage. As we have just seen, adding knowledge of $W_T$ does
not destroy the semimartingale property. But in the Black--Scholes
model (with constant coefficients), knowing $W_T$ is the same as
knowing the terminal stock price $S_T$, and that clearly entails
arbitrage.
\endbox
\end{remarks}

Consider a probability space Let $(\Omega,{\cal F},\mathbb P)$
equipped with right--continuous filtrations $\mathbb F = ({\cal
F}_t)_{t\geq 0},\mathbb H = ({\cal H}_t)_{t\geq 0}$. We are
interested in the enlarged filtration $\mathbb G = ({\cal
G}_t)_{t\geq 0}$ which amalgamates the information in $\mathbb
F,\mathbb H$. The natural candidate for this is $({\cal
F}_t\lor{\cal H}_t)_t := (\sigma({\cal F}_t\cup{\cal H}_t)_t$.
However, in order to do analysis, we need to work with c\`adla\`ag
versions of stochastic processes, and for those to exist, we require
the usual hypotheses to hold. Hence we will define
\[{\cal G}_t = \bigcap_{s>t} ({\cal F}_s\lor{\cal H}_s)\]
Two special cases of enlargement have received significant
attention:
\begin{enumerate}[(1)]\item If $X$ is a random variable, and ${\cal
H}_t = \sigma(X)$ for all $t$, then \[{\cal G}_t=\bigcap_{s>t}({\cal
F}_s\lor\sigma(X))\] is an {\em initial enlargement} of $\mathbb F$
by $X$. The information $\sigma(X)$ is added to $\mathbb F$ all at
once.
\item If $T$ is a random time, but not necessarily an $\mathbb F$--stopping
time, then we can gradually add just enough information to make it a
stopping time. Then \[{\cal G}_t = \bigcap_{s>t}({\cal
F}_s\lor\sigma(T\land s))\] Here $\mathbb G$ is a {\em progressive
enlargement} of $\mathbb F$.\newline(For example, an insider may not
know exactly when an event (e.g. bankruptcy) is going to happen, but
will know if it has happened. The time $T$ of bankruptcy is then a
stopping time for the insider, but not the \lq\lq honest" investor.
\end{enumerate}

In this chapter, we are mainly concerned with initial enlargements.
nevertheless, we start off in a very general framework.

\section{Enlargements of Filtrations via Changes of Measure}
Let $(\Omega,{\cal F},\mathbb P)$ be a probability space with
right--continuous filtrations $\mathbb F = ({\cal F}_t)_{t\geq
0},\mathbb H = ({\cal H}_t)_{t\geq 0}$. The object of our studies is
the enlarged filtration $\mathbb G = ({\cal G}_t)_{t\geq 0}$, where
\[{\cal G}_ t = \bigcap_{s>t}({\cal F}_s\lor {\cal H}_s)\] Following
Ankirchner et al.(\cite{ADI2005},\cite{Ankirchner2005}), we find a
\lq\lq translation" of the current set--up in a product space, where
it will transpire that the enlargement of filtration corresponds to
a change of measure. Consider, therefore, the product space
$(\bar{\Omega}, \bar{\cal F},\bar{\mathbb P},\bar{\mathbb F})$,
where
\[\bar{\Omega} = \Omega\times\Omega\qquad \bar{\cal F} = {\cal
F}_\infty\otimes{\cal H}_\infty\qquad\bar{\mathbb P} = \mathbb
Pi^{-1}\qquad \bar{\cal F}_t = \bigcap_{s>t}({\cal
F}_s\otimes{\cal H}_s)\] and $
i:\Omega\to\bar{\Omega}:\omega\mapsto(\omega,\omega)$ is the
diagonal embedding. By the change of variable formula, it follows
that\[\int f(\omega,\omega')\;d\bar{\mathbb P}(\omega,\omega') =
\int f(\omega,\omega)\;d\mathbb P(\omega)\] for every $\bar{\cal
F}$--measurable $f$.

Right now,  our main concern is how various properties from
stochastic analysis fare under the translation. We have two
translation mappings, the diagonal embedding, and the first
projection mapping
\[ i:\Omega\to\bar{\Omega}:\omega\mapsto(\omega,\omega)\qquad
\pi:\bar{\Omega}\to\Omega:(\omega,\omega')\mapsto\omega\] The map $
i$ is used to translate from $\bar{\Omega}$ to $\Omega$. We will
show, roughly, that if $\bar{X}$ is
measurable/predictable/semimartingale w.r.t $(\bar{\mathbb F},\bar{
\mathbb P})$ then $\bar{X}\circ  i$ is
measurable/predictable/semimartingale w.r.t $(\mathbb G,\mathbb P)$.
Conversely, the map $\pi$ is used to translate from $\Omega$ to
$\bar{\Omega}$, and we will see, roughly, that if $X$ is a $(\mathbb
F,\mathbb P)$--semimartingale, then $X\circ\pi$ is a $(\bar{\mathbb
F},\bar{\mathbb P})$--semimartingale. We say \lq\lq roughly",
because we do have to keep track of null sets. We will be working
with different measures, which therefore have different null sets.
Now much of stochastic analysis works only when we impose the usual
conditions on filtrations --- e.g., these are required to prove that
martingales have c\`adl\`ag versions --- so that we must ensure that
our filtrations are properly augmented w.r.t. the right measure, if
we want to be able to do any analysis. Thus, given a filtration
$\mathbb K=({\cal K}_t)_{t\geq 0}$ and a probability measure $P$, we
let $\mathbb K^P = ({\cal K}^P_t)_{t\geq 0}$ be the filtration
$\mathbb K$ completed by $P$--negligible sets.

\begin{lemma}\label{lemma_Omegabar_to_Omega_G}\begin{enumerate}[(a)]\item
We have
\[{\cal G}_t = i^{-1}(\bar{\cal F}_t)\qquad {\cal G}_t^\mathbb P
\supseteq i^{-1}(\bar{\cal F}_t^{\bar{\mathbb P}})\] \item If
$\bar{X}_n\to 0$ in $\bar{\mathbb P}$--probability, then
$\bar{X}_n\circ i\to 0$ in $\mathbb P$--probability. \item If
$\bar\tau$ is a $\bar{\mathbb F}^{\bar{\mathbb P}}$--stopping
time, then $\bar{\tau}\circ i$ is a $\mathbb G^\mathbb
P$--stopping time.
\end{enumerate}\end{lemma}

\bproof (cf. \cite{Ankirchner2005})  (a) As ${\cal F}_t\lor{\cal
H}_t =\sigma(F\cap H: F\in{\cal F}_t, H\in{\cal H}_t)$, we have
\[\aligned {\cal G}_t &=\bigcap_{s>t}\sigma(F\cap H: F\in{\cal F}_s, H\in{\cal
H}_s)\\&=\bigcap_{s>t}\sigma(i^{-1}(F\times H):F\in{\cal F}_s,
H\in{\cal H}_s)\\&=i^{-1}\left(\bigcap_{s>t}{\cal F}_s\otimes{\cal
H}_s\right)\\&=i^{-1}(\bar{\cal F}_t)
\endaligned\]
Now if $A\in\bar{\cal F}_t^{\bar{\mathbb P}}$, then there exists
$B\in\bar{\cal F}_t$ such that $\bar{\mathbb P}(A\Delta B) = 0$.
Hence $\mathbb P(i^{-1}(A)\Delta i^{-1}(B))=0$. Now
$i^{-1}(B)\in{\cal G}_t$, and hence $i^{-1}(A)\in{\cal
G}_t^\mathbb P$. \newline (b) Let $X_n= \bar{X_n}\circ i$. Then
$\mathbb P(|X_n|>\varepsilon) = \mathbb
P(i^{-1}\{|\bar{X}_n|>\varepsilon\}) = \bar{\mathbb
P}(|\bar{X}_n|>\varepsilon)\to 0$ as $n\to\infty$.\newline (c) If
$\bar\tau$ is a $\bar{\mathbb F}^{\bar{\mathbb P}}$--stopping
time, then
\[\{\bar{\tau}\circ i\leq t\} = i^{-1}\{\bar\tau\leq t\}\in
i^{-1}(\bar{\cal F}^{\bar{\mathbb P}}_t)\subseteq {\cal
G}^{\mathbb P}_t\] so that $\bar\tau\circ i$ is a $\mathbb
G^\mathbb P$--stopping time.\eproof

 The following result translates stochastic analysis from $(\bar{\mathbb F}^{\bar{\mathbb P}},\bar{\mathbb P})$
to $(\mathbb G^\mathbb P,\mathbb P)$:

\begin{theorem}\label{thm_Omegabar_to_Omega_G}\begin{enumerate}[(a)]\item If
$\bar{\Omega}\buildrel{\bar{f}}\over\to \mathbb R$ is
$\bar{\mathbb F}^{\bar{\mathbb P}}$--measurable, then
$\bar{f}\circ i$ is $\mathbb G^\mathbb P$--measurable. \item If
$\bar{X}$ is $\bar{\mathbb F}^{\bar{\mathbb P}}$--adapted, then
$\bar{X}\circ i$ is $\mathbb G^\mathbb P$--adapted. \item If
$\bar{X}$ is $\bar{\mathbb F}^{\bar{\mathbb P}}$--predictable,
then $\bar{X}\circ i$ is $\mathbb G^\mathbb P$--predictable. \item
If $\bar{X}$ is a $(\bar{\mathbb F}^{\bar{\mathbb P}},\bar{\mathbb
P})$--(local) martingale, then $\bar{X}\circ i$ is a $(\mathbb
G^\mathbb P,\mathbb P)$--(local) martingale. \item If $\bar{X}$ is
a $(\bar{\mathbb F}^{\bar{\mathbb P}},\bar{\mathbb
P})$--semimartingale, then $\bar{X}\circ i$ is a $(\mathbb
G^\mathbb P,\mathbb P)$--semimartingale.

\item If $\bar{H}$ is c\`agl\`ad $(\bar{\mathbb F}^{\bar{\mathbb
P}},\bar{\mathbb P})$--adapted and $\bar{X}$ a $(\bar{\mathbb
F}^{\bar{\mathbb P}},\bar{\mathbb P})$--semimartingale, and $H =
\bar{H}\circ i, X=\bar{X}\circ i$, then
\[(\bar{H}\bullet\bar{X})\circ i = H\bullet X\] up to
indistinguishability. \item If $\bar{X},\bar{Y}$ are
$(\bar{\mathbb F}^{\bar{\mathbb P}},\bar{\mathbb
P})$--semimartingales, and $X=\bar{X}\circ i, Y = \bar{Y}\circ i$,
then
\[[\bar{X},\bar{Y}]\circ i = [X,Y]\]
up to indistinguishability.
\end{enumerate}
\end{theorem}

\bproof (a) follows from Lemma \ref{lemma_Omegabar_to_Omega_G} by
the usual arguments, and (b) follows immediately from (a).\newline
(c) Consider first an $(\bar{\mathbb F}^{\bar{\mathbb
P}},\bar{\mathbb P})$--predictable process of the form $\bar{X} =
I_AI_{(s,t]}$, where $A\in\bar{\cal F}^{\bar{\mathbb P}}_s$. Then
$I_A\circ i$ is ${\cal G}^{\mathbb P}_s$--measurable, by (a), and
so $\bar{X}\circ i = (I_A\circ i)I_{(s,t]}$ is $(\mathbb
G^{\mathbb P},\mathbb P)$--predictable. The result for general
predictable $\bar{X}$ now follows using a monotone class
theorem.\newline (d)
 Suppose first that $\bar{X}$ is a $(\bar{\mathbb F}^{\bar{\mathbb P}},\bar{\mathbb
P})$--martingale. Put $X = \bar{X}\circ i$ and let $0\leq s<t,
G\in{\cal G}^{\mathbb P}_s$. To show that $X$ is a $(\mathbb
G^{\mathbb P},\mathbb P)$--martingale, it suffices to show that
$\mathbb E_{\mathbb P}[I_G(X_t-X_s)]=0$. Now by Lemma
\ref{lemma_Omegabar_to_Omega_G}, there is $\bar{F}\in\bar{\cal
F}^{\bar{\mathbb P}}_s$ such that $G = i^{-1}(\bar{F})$. Then
\[\mathbb E_{\mathbb P}[I_G(X_s-X_t)] =
\int (I_{\bar{F}}(\bar{X}_t-\bar{X}_s))\circ i\;d\mathbb P = \int
I_{\bar{F}}(\bar{X}_t-\bar{X}_s)\;d\bar{\mathbb P} =\mathbb
E_{\bar{\mathbb P}}[I_{\bar{F}}(\bar{X}_t-\bar{X}_s)]=0\] If
$\bar{X}$ is a $(\mathbb G^{\mathbb P},\mathbb P)$--local
martingale, the result follows easily by localization, using (d).
\newline (e) We sketch the proof. Suppose that $\bar{X}$ is a $(\bar{\mathbb F}^{\bar{\mathbb P}},\bar{\mathbb
P})$--semimartingale, and put $X = \bar{X}\circ i$. Clearly, $X$
is c\`adl\`ag, because $\bar{X}$ is, and $\mathbb G^{\mathbb
P}$--adapted. To prove that $X$ is a $(\mathbb G^\mathbb P,\mathbb
P)$--semimartingale, it suffices to show that if $(\theta^n)_n$ is
a sequence of simple predictable processes converging to zero
uniformly in $(t,\omega)$, then $(\theta^n\bullet Y)$ converges to
zero in $\mathbb P$--probability, by the Bichteler--Dellacherie
Theorem.  Now if $\theta =\sum_{i=1}^n\theta_iI_{(t_i,t_{i+1}]}$
is simple, where each $\theta_i$ is ${\cal G}_{t_i}$--measurable,
then $\theta_i = \bar{\theta_i}\circ i$ for some $\bar{\cal
F}_{t_i}$--measurable $\bar{\theta}_i$. Consequently, for every
simple $\mathbb G$--predictable sequence $(\theta^n)_n$ converging
to 0 uniformly, we can find a simple $\bar{\mathbb
F}$--predictable sequence $(\bar{\theta}^n)_n$ converging to 0
uniformly such that $\theta^n = \bar{\theta^n}\circ i$: Since
$\bar X$ is a semimartingale, $(\bar{\theta}^n\bullet\bar{X})$
converges to 0 in $\bar{\mathbb P}$--probability, and hence
$(\theta^n\bullet X)$ converges to 0 in $\mathbb P$--probability.
\newline (f) For $n\in\mathbb N$, and $1\leq i\leq 2^n$, let
$t^n_i = ti2^{-n}$. Then
\[\bar{H}_0\bar{X}_0 +\sum_{i=0}^{2^n-1} \bar{H}_{t^n_i}(\bar{X}_{t^n_{i+1}}-\bar{X}_{t^n_i})\]
converges to $(\bar{H}\bullet\bar{X})_t$ in $\bar{\mathbb
P}$--probability, as $n\to\infty$. Thus
$(\bar{H}\bullet\bar{X})_t\circ\ i$ is the limit in $\mathbb
P$--probability of
\[{H}_0{X}_0 +\sum_{i=0}^{2^n-1}
{H}_{t^n_i}({X}_{t^n_{i+1}}-{X}_{t^n_i})\] This limit is also
$(H\bullet X)_t$, and so $(\bar{H}\bullet\bar{X})\circ\ i$ and
$(H\bullet X)$ are modifications of each other. Since  both are
c\`adl\`ag, they must be indistinguishable under $\mathbb
P$.\newline (g) follows in a similar way, using the fact that
$[\bar{X},\bar{Y}]_t$ is a limit in probability of
$\bar{X}_0\bar{Y}_0 +\sum_{i=0}^{2^n-1}
(\bar{X}_{t^n_{i+1}}-\bar{X}_{t^n_i})(\bar{Y}_{t^n_{i+1}}-\bar{Y}_{t^n_i})$
 \eproof

To translate stochastic analysis from $\Omega$ to $\bar{\Omega}$,
we don't need quite as much. Given a probability measure $R$ on
$(\Omega,{\cal F})$, let
\[\bar{\mathbb Q}=\mathbb P|_{{\cal F}_\infty}\otimes R|_{{\cal
H}_\infty}\] $\bar{\mathbb Q}$ is called a {\em decoupling}
measure. We shall often take $R=\mathbb P$.

\begin{lemma} If $A\in{\cal F}^{\mathbb P}_t$, then
$A\times \Omega\in\bar{\cal F}^{\bar{\mathbb Q}}_t$.
\end{lemma}
\bproof Take $B$ in ${\cal F}_t$ so that $\mathbb P(A\Delta B) =
0$, and note that $\bar{\mathbb Q}((A\times \Omega)\Delta
(B\times\Omega)) = \mathbb P(A\Delta B) = 0$.\eproof

\begin{theorem}\label{thm_M_to_Mbar_local_mg} Let $M$ be a right--continuous (everywhere) $(\mathbb F^\mathbb P,\mathbb
P)$--(local) martingale. Then $M\circ\pi$ is a $(\bar{\mathbb
F}^{\bar{\mathbb Q}},\bar{\mathbb Q})$--(local) martingale.
\end{theorem}

\bproof Let $\bar{M}=M\circ \pi$, i.e. $\bar{M}(\omega,\omega') =
M(\omega)$. It is clear that $\bar{M}$ is right--continuous,
because $M$ is. To see that $\bar{M}$ is $\bar{\mathbb
F}^{\bar{\mathbb Q}}$--adapted,  note that $c\in\mathbb R$, we
have
\[\{\bar{M}_t\leq c\}=\{M_t\leq c\}\times\Omega \in\bar{\cal F}^{\bar{\mathbb
Q}}_t\] by the lemma, which proves that $\bar{M}$ is adapted.

Suppose first that $M$ is a martingale. For $0\leq s<t$,
$A\in{\cal F}_s, B\in{\cal H}_s$, we have \[\mathbb
E_{\bar{\mathbb Q}}[I_{A\times B}(\bar{M}_t-\bar{M}_s)] =
R(B)\mathbb E_{\mathbb P}[I_A(M_t-M_s)] = 0\] By a monotone class
argument, it follows that $\mathbb E_{\bar{\mathbb
Q}}[\bar{X}(\bar{M}_t-\bar{M}_s)] = 0$ for all bounded ${\cal
F}_s\otimes{\cal H}_s$--measurable $\bar{X}$. As $\bar{M}$ is
right--continuous, we see that $\mathbb E_{\bar{\mathbb
Q}}[\bar{X}(\bar{M}_t-\bar{M}_s)] = 0$ for all bounded
$\bigcap_{u>s}({\cal F}_u\otimes{\cal H}_u)$--measurable
$\bar{X}$, which proves that $\bar{M}$ is a $(\bar{\mathbb
F}^{\bar{\mathbb Q}},\bar{\mathbb Q})$--martingale.

The result can now easily be extended (via localization) to the
case where $M$ is a local martingale: If $\tau$ is a $\mathbb
F^\mathbb P$--stopping time, and $\bar{\tau}=\tau\circ \pi$, then
\[\{\bar{\tau}\leq t\} = \{\tau\leq t\}\times\Omega\in\bar{\cal
F}^{\bar{\mathbb Q}}_t\] by the lemma, so that $\bar{\tau}$ is a
$\bar{\mathbb F}^{\bar{\mathbb Q}}$--stopping time.\eproof

We shall henceforth impose the following assumption on
$\bar{\mathbb P},\bar{\mathbb Q}$:\begin{center}\fbox{{\bf
Assumption $\mathbf A$:}\qquad\qquad $\bar{\mathbb
P}<<\bar{\mathbb Q}$  \qquad on\qquad $\bar{\cal
F}$\phantom{aaaaaaaaaaaaaaaaaaaaaaaaa}}\end{center} Now recall
that the semimartingale property is preserved under a change of
measure (provided the new measure is absolutely continuous w.r.t.
the original). In particular, every $(\bar{\mathbb
F}^{\bar{\mathbb Q}},\bar{\mathbb Q})$--semimartingale is a
$(\bar{\mathbb F}^{\bar{\mathbb Q}},\bar{\mathbb
P})$--semimartingale. We therefore immediately obtain the
following important result:

 \begin{theorem}\label{thm_A_implies H'} $(\mathbf A)$ implies $(\mathbf H')$. \newline I.e. if
 $\bar{\mathbb P}<<\bar{\mathbb Q}$, then every $(\mathbb F^{\mathbb
 P},\mathbb P)$--semimartingale is a $(\mathbb G^{\mathbb P},\mathbb
 P)$--semimartingale.
 \end{theorem}

\bproof Suppose that $M$ is a $(\mathbb F^{\mathbb
 P},\mathbb P)$--local martingale. It suffices to show that it is also a $(\mathbb G^{\mathbb P},\mathbb
 P)$--semimartingale. Now $\bar{M}=M\circ\pi$ is a $(\bar{\mathbb F}^{\bar{\mathbb
Q}},\bar{\mathbb Q})$--local martingale, by Theorem
\ref{thm_M_to_Mbar_local_mg}. Since $\bar{\mathbb P}<<\bar{\mathbb
Q}$, we have that $\bar{M}$ is a $(\bar{\mathbb F}^{\bar{\mathbb
Q}},\bar{\mathbb P})$--semimartingale. By Stricker's Theorem, it
is also a $(\bar{\mathbb F}^{\bar{\mathbb P}},\bar{\mathbb
P})$--semimartingale (because $\bar{M}$ is adapted to
$\bar{\mathbb F}^{\bar{\mathbb P}}\subseteq \bar{\mathbb
F}^{\bar{\mathbb Q}}$). Hence $\bar{M}\circ i = M\circ\pi\circ i =
M$ is a $(\mathbb G^\mathbb P,\mathbb P)$--semimartingale, by
Theorem \ref{thm_Omegabar_to_Omega_G}.\eproof

\section{Doob--Meyer Decompositions via Girsanov's Theorem} We have
now accomplished the following: Given a probability space
$(\Omega,{\cal F},\mathbb P)$  with filtrations $\mathbb F =
({\cal F}_t)_{t\geq 0},\mathbb H = ({\cal H}_t)_{t\geq 0}$
satisfying the usual conditions, we have shown that assumption
$(\mathbf A)$ implies hypothesis $(\mathbf H')$ for the enlarged
filtration $\mathbb G = (\bigcap_{s>t}({\cal F}_s\lor {\cal
H}_s))_{t\geq 0}$. In essence, ignoring negligible sets, we
embedded the original set--up $(\Omega,{\cal F},\mathbb P,\mathbb
F)$ into a product space
\[(\bar{\Omega},\bar{\cal F},\bar{\mathbb Q}, \bar{\mathbb F}) =
(\Omega\times\Omega, {\cal F}_\infty\otimes{\cal H}_\infty,
\mathbb P|{\cal F}_\infty\otimes\mathbb P|{\cal H}_\infty,
(\bigcap_{s>t}{\cal F}_s\otimes{\cal H}_s)_{t\geq 0})\] The first
projection $\pi:(\bar{\Omega},\bar{\cal F}_t,\bar{\mathbb
Q})\to(\Omega,{\cal F}_t,\mathbb P)
:(\omega,\omega')\mapsto\omega$ is clearly measurable (for each
$t\geq 0$), i.e. the first component captures the initial set--up
$(\Omega,{\cal F},\mathbb P,\mathbb F)$. Theorem
\ref{thm_Omegabar_to_Omega_G} however, shows that the enlarged
set--up $(\Omega,{\cal F},\mathbb P,\mathbb G)$ is captured
perfectly by the space $(\bar{\Omega},\bar{\cal F},\bar{\mathbb
P}, \bar{\mathbb F})$. In essence, therefore, enlarging
filtrations on $\Omega$ corresponds to changing measures on
$\bar{\Omega}$. This means that we can bring the machinery of
Girsanov's Theorem into play. Again, we rely entirely on
\cite{Ankirchner2005},\cite{ADI2005}.

Since we are assuming that $\bar{\mathbb P}<<\bar{\mathbb Q}$ on
$\bar{\cal F}$, let $\bar{Z}$ be a c\`adl\`ag version of the
likelihood process:
\[\bar{Z_t} = \frac{d\bar{\mathbb P}}{d\bar{\mathbb
Q}}\Big|\bar{\cal F}^{\bar{\mathbb Q}}_t\] (The regularization
theorems guaranteeing the existence of a c\`adl\`ag version
require a right--continuous and complete filtration.) As $\bar{Z}$
is a $(\bar{\mathbb F}^{\bar{\mathbb Q}},\bar{\mathbb
Q})$--martingale, it follows that the process $Z=\bar{Z}\circ i$
is a $(\mathbb G^\mathbb P,\mathbb P)$--semimartingale.

Now the Lenglart--Girsanov theorem states that every $\bar{\mathbb
F}^{\bar{\mathbb Q}},\bar{\mathbb Q})$--semimartingale is also a
$\bar{\mathbb F}^{\bar{\mathbb Q}},\bar{\mathbb
P})$--semimartingale. To be precise, let $\bar{M}$ be a
$(\bar{\mathbb F}^{\bar{\mathbb Q}},\bar{\mathbb Q})$--local
martingale with $\bar{M}_0 = 0$. Define
\[\bar{T} = \inf\{t>0:\bar{Z}_t = 0, \bar{Z}_{t-}>0\}\qquad \bar{U}_t =\Delta\bar{M}_{\bar{T}}I_{\{t\geq\bar{T}\}}\]
and let $\tilde{U}$ be the $(\bar{\mathbb F}^{\bar{\mathbb
Q}},\bar{\mathbb Q})$--compensator of $\bar{U}$ (i.e. the {\em
unique predictable} finite variation process such that $\bar{U}
-\tilde{U}$ is a $(\bar{\mathbb F}^{\bar{\mathbb Q}},\bar{\mathbb
Q})$--local martingale; cf. Protter\cite{Protter2004}, III.5). By
the Lenglart--Girsanov Theorem (cf. Protter\cite{Protter2004},
III.8), the process
\[\bar{M}_t - \int_0^t\frac{1}{\bar{Z}_s}\;d[\bar{M},\bar{Z}]_s
+\tilde{U_s}\] is a $(\bar{\mathbb F}^{\bar{\mathbb P}},\bar{\mathbb
P})$--local martingale. Now let $U= \tilde{U}\circ i$, and recall
that $Z=\bar{Z}\circ i$. Then:

\begin{theorem} If $M$ is a $(\mathbb F^{\mathbb P},\mathbb
P)$--local martingale with $M_0=0$, then
\[M-\frac{1}{Z}\bullet[Z,M]+U\]
is a $(\mathbb G^{\mathbb P},\mathbb P)$--local martingale.
\end{theorem}

\bproof  By Theorem \ref{thm_M_to_Mbar_local_mg}, the process
$\bar{M}=M\circ\pi$ is a $(\bar{\mathbb F}^{\bar{\mathbb
Q}},\bar{\mathbb Q})$--local martingale with $\bar{M}_0 = 0$. By the
Lenglart--Girsanov Theorem,
\[\bar{X}=
\bar{M}-\frac{1}{\bar{Z}}\bullet[\bar{Z},\bar{M}]+\tilde{U}\] is a
$(\bar{\mathbb F}^{\bar{\mathbb P}},\bar{\mathbb P})$--local
martingale. By Theorem \ref{thm_Omegabar_to_Omega_G}, $\bar{X}\circ
i$ is a $(\mathbb G^\mathbb P,\mathbb P)$--local martingale. But, by
the same theorem, we see that
\[X\circ i = M-\frac{1}{Z}\bullet[Z,M]+U\]
\eproof

If $M$ is a continuous $(\mathbb F^{\mathbb P},\mathbb P)$--local
martingale with $M_0=0$, and $\bar{M} = M\circ\pi$, etc., then
clearly $U = 0$ (because $\bar{U} =0$, and hence $\tilde{U}=0$), so
that we have:
\begin{theorem} If $M$ is a continuous $(\mathbb F^{\mathbb P},\mathbb
P)$--local martingale with $M_0=0$, then
\[M-\frac{1}{Z}\bullet[Z,M]\]
is a $(\mathbb G^{\mathbb P},\mathbb P)$--local martingale.
\end{theorem}

\section{Jacod's Criterion for Initial Enlargements}
Let $(\Omega,{\cal F},\mathbb P,\mathbb F)$ be a filtered
probability space satisfying the usual hypotheses, and let $X$ be an
${\cal F}$--measurable random element with values in a state space
$(S,{\cal S})$. We consider an initial enlargement $\mathbb G$ of
$\mathbb F$, given by
\[{\cal G}_t=\bigcap_{s>t}({\cal F}_s\lor\sigma(X))\]
i.e. \[{\cal G}_t=\bigcap_{s>t}({\cal F}_s\lor{\cal
H}_s)\qquad\text{where }\quad{\cal H}_s = \sigma(X)\quad\text{ for
all }s\geq 0\] The most well--known condition for an initial
enlargement to satisfy condition $(\mathbf H')$ is {\em Jacod's
Criterion} (due to Jacod\cite{Jacod1985}), which we shall now prove.

First recall that if $X$ is a random element with values in a Borel
space $(S,{\cal S})$, then $X$ has a (unique) regular conditional
distribution $Q_t(\omega,dx)$ with the properties
that\begin{enumerate}[(i)]\item $A\mapsto Q_t(\omega, A)$ is a
probability measure on $(S,{\cal S})$, for almost all $\omega$.
\item $\omega\mapsto Q_t(\omega, A)$ is measurable for all
$A\in{\cal S}$;
\item $Q_t(\cdot, dx)$ is a version of $\mathbb P(X\in dx|{\cal
F}_t)$.
\end{enumerate}

Jacod's criterion depends on the following assumption:
\begin{center}\fbox{\begin{minipage}{12cm}{\bf
Assumption $\mathbf J$:} For each $t\geq 0$, here exists a
$\sigma$--finite measure $\eta_t(ds)$ on the state space $(S,{\cal
S})$ of $X$ such that
\[Q_t(\omega,ds)<<\eta_t(ds)\quad\mathbb P\text{--a.s.}\]\end{minipage}}\end{center}
\begin{theorem}\label{thm_Jacod_criterion} {\rm (Jacod's Criterion}) Let $X$ be a random element on $(\Omega,{\cal F},\mathbb P,\mathbb F)$ with regular conditional
distribution $Q_t(\omega,dx)$. If $(\mathbf J)$ holds, then
$(\mathbf H')$ holds for the initial enlargement of $\mathbb F$ with
$\sigma(X)$.
\end{theorem}

We prove this result in two steps. First we show that $(\mathbf J)$
is equivalent to a seemingly stronger statement $(\mathbf J')$:
\begin{center}\fbox{\begin{minipage}{12cm}{\bf
Assumption $\mathbf J'$:} there exists a {\em single}
$\sigma$--finite measure $\eta(ds)$ on the state space $(S,{\cal
S})$ of $X$ such that
\[Q_t(\omega, ds)<<\eta(ds)\qquad\text{ for all
}\omega\in\Omega\text{ and }t\geq 0\]\end{minipage}}\end{center}
Then we show that $(\mathbf J')$ is equivalent to $(\mathbf A)$.

\begin{lemma} \label{lemma_abs_continuity_law_X} In the set--up of Theorem
\ref{thm_Jacod_criterion}, if $(\mathbf J)$ holds, then so does
$(\mathbf J')$.\newline In fact we may take the {\em single}
$\sigma$--finite measure $\eta$ on the state space $(S,{\cal S})$ of
$X$ satisfying
\[Q_t(\omega, dx)<<\eta(dx)\qquad\text{ for all
}\omega\in\Omega\text{ and }t\geq 0\] to be the law of $X$.
\end{lemma}

\bproof We follow Protter\cite{Protter2004}. Assume that $(\mathbf
J')$ holds,, and that $Q_t(\omega,ds)<<\eta_t(ds)$ for all $\omega$.
By Theorem \ref{thm_Doob_disintegration} there exists a ${\cal
F}_t\otimes{\cal S}$--measurable function $q_t(\omega,s)$ such that
$Q_t(\omega,ds) = q_t(\omega, s)\;\eta(ds)$. Define $a_t(s) =
\mathbb E[q_t(\cdot,s)]$ and define
\[r_t(\omega,s) = \left\{\aligned \frac{q_t(\omega,s)}{a_t(s)}\qquad
&\text{if }a_t(s)>0\\0\qquad\qquad&\text{else}\endaligned\right.\]
Since we must have $q_t(\omega,s) = 0$ a.s. whenever $a_t(s) = 0$,
it follows that $q_t(\omega,s) = r_t(\omega,s)a_t(s)$ a.s. Hence
\[Q_t(\omega,ds) = r_t(\omega,s)a_t(s)\eta_t(ds)\qquad\text{a.s.}\]
Now let $\eta$ be the law of $X$. For every non--negative ${\cal
S}$--measurable function $g$, and every $t\geq 0$, we have
\[\aligned \int g(s)\;\eta(ds)&= \mathbb E[g(X)]\\&=\mathbb
E\left[\int g(s)\; Q_t(\cdot,ds)\right]\\&=\mathbb E\left[\int
g(s)q_t(\cdot,s)\eta_t(ds)\right]\\&=\int g(s)\mathbb
E[q_t(\cdot,s)]\;\eta_t(ds)\\&=\int
g(s)a_t(s)\;\eta_t(ds)\endaligned\] Hence\[a_t(s)\eta_t(ds) =
\eta(ds)\qquad\text{for all }t\geq 0\] so that
$Q_t(\omega,ds)=r_t(\omega,s)\eta(ds)$, as required.\eproof

\begin{lemma}\label{lemma_A is_J'_1}\begin{enumerate}[(i)]\item
The map $j:(\Omega\times\Omega,{\cal F}_t\otimes \sigma(X))\to
(\Omega\times S,{\cal F}_t\otimes{\cal S}):(\omega,\omega')\mapsto
(\omega, X(\omega'))$ is measurable. \item Every $C\in{\cal
F}_t\otimes\sigma(X)$ is of the form $j^{-1}(D)$ for some $D\in
{\cal F}_t\otimes{\cal S}$.
\end{enumerate}\end{lemma}
\bproof (i) follows from the fact that $j^{-1}(A\times B)=A\times
X^{-1}(B)$. \newline (ii) The family of all $C$ which can be
represented as $j^{-1}(D)$, is a $\sigma$--algebra containing all
the measurable rectangles in ${\cal F}_t\otimes \sigma(X)$.\eproof

\begin{lemma}\label{lemma_A_is_J'_2}
If $f:(\Omega\times S,{\cal F}_t\otimes{\cal S})\to\mathbb R$ is
measurable, then
\[\int_{\Omega\times\Omega} f(\omega, X(\omega'))\;\bar{\mathbb
P}(d\omega) = \int_\Omega\int_S f(\omega,s)Q_t(\omega,ds)\;\mathbb
P(d\omega)\]
\end{lemma}
\bproof Define a probability measure $\tilde{P}$ on ${\cal
F}_t\otimes \sigma(X)$ by
\[\tilde{P}(A\times X^{-1}(B)) =\int_A Q_t(\omega,B)\;\mathbb P(d\omega)\qquad\text{for }A\in{\cal F}_t, B\in {\cal S}\]\eproof
On such rectangles $A\times X^{-1}(B)$,  we have
\[\tilde{P}(A\times X^{-1}(B)) = \mathbb E[Q_t(\cdot, B)I_A] =
\mathbb E[I_{X\in B}I_A]=\mathbb P(A\cap X^{-1}(B)) = \bar{\mathbb
P}(A\times X^{-1}(B))\] Consequently, $\tilde{P}=\bar{\mathbb P}$.

It is now easy to verify that
\[\int f(\omega,X(\omega'))\;\tilde{P}(d\omega\times d\omega')
=\iint f(\omega,s)\;Q_t(\omega,ds)\;\mathbb P(d\omega)\] This
follows directly from the definition of $\tilde{P}$ if $f=I_{A\times
B}$ is the indicator of a measurable rectangle, and the follows for
arbitrary measurable $f$ by the usual arguments. \eproof

\begin{proposition}\label{propn_A_is_J'} In the set--up of Theorem \ref{thm_Jacod_criterion},
$(\mathbf J')$ is equivalent to $(\mathbf A)$: To be precise, let
$R$ be a probability measure with $RX^{-1}=\eta$, where $\eta$ is
the emasure supplied by $(\mathbf J')$. (In particular, we may take
$R=\mathbb
 P$, because we may take $\eta$ to be $\mathbb PX^{-1}$, by Lemma \ref{lemma_abs_continuity_law_X}.)
 \newline With ${\cal H}_t=\sigma(X)$, we have
\[\bar{\mathbb P}<<\bar{\mathbb Q}:=\mathbb P\otimes\mathbb R\quad\text{on
}\bar{\cal F}_t\text{ for all }t\geq 0\]if and only
if\[Q_t(\omega,ds)<<\eta(ds)\quad\text{for almost all $\omega$, for
all $t\geq 0$}\]
\end{proposition}
\bproof (\cite{Ankirchner2005}, \cite{ADI2005}) First assume
$(\mathbf J')$, and thus that the regular conditional version
$Q_t(\omega,ds)$ of $\mathbb P(X\in ds|{\cal F}_t)$ is absolutely
continuous w.r.t. the law $\eta$ of $X$, for each $t\geq 0$ --- cf.
Lemma \ref{lemma_abs_continuity_law_X}. Let $t\geq 0$.

Now, for $s<t$, let $C\in\bar{\cal F}_s=\bigcap_{u>s}({\cal
F}_u\otimes\sigma(X))$ with $\bar{\mathbb Q}(C) = 0$. By Lemma
\ref{lemma_A is_J'_1}, there is $D\in{\cal F}_t\otimes{\cal S}$ such
that $C=j^{-1}(D)$. Then, since $\bar{\mathbb Q}=\mathbb
P\otimes\mathbb P$ (restricted to suitable $\sigma$--algebras)
\[\aligned \bar{\mathbb Q}(C) &= \int I_D\circ
j(\omega,\omega')\;\bar{\mathbb Q}(d\omega\times d\omega')\\
&= \int_\Omega\left(\int_SI_D(\omega,s)\;\eta(ds)\right)\;\mathbb
P(d\omega)=0
\endaligned\]
The inner integral (being non--negative) must therefore satisfy
$\int_SI_D(\omega,s)\;\eta(ds) = 0$ for $\mathbb P$--a.a. $\omega$.
As $Q_t(\omega, ds)<<\eta(ds)$, also
$\int_SI_D(\omega,s)\;Q_t(\omega,ds) = 0$ for almost all $t,
\omega$. Hence by Lemma \ref{lemma_A_is_J'_2},
\[\bar{\mathbb P}(C) =  \int_\Omega I_D(\omega,X(\omega'))\;\bar{\mathbb P}(d\omega\times d\omega')
=\int_\Omega\int_S I_D(\omega,s)Q_t(\omega,ds)\;\mathbb P(d\omega) =
0\] This proves that $(\mathbf J')\Rightarrow(\mathbf A)$.

Conversely, suppose that $\bar{\mathbb P}<<\bar{\mathbb Q} = \mathbb
P\otimes R$ on $\bar{\cal F}_t$.  If $\varphi(\omega,\omega')$ is a
${\cal F}_t\otimes\sigma(X)$--measurable version of
$\frac{d\bar{\mathbb P}}{d\bar{\mathbb Q}}$, then by lemma
\ref{lemma_A is_J'_1}, there is a ${\cal F}_t\otimes{\cal
S}$--measurable $\tilde{\varphi}$ such that $\varphi(\omega,
\omega') = \tilde{\varphi}(\omega, X(\omega'))$. Define $U_t(\omega,
B)=\int_{X^{-1}(B)}\tilde\varphi(\omega,X(\omega'))\;R(d\omega')$.
It is easy to see that $U_t$ is  a stochastic kernel from
$(\Omega,{\cal F}_t)$ to $(S,{\cal S})$. We claim that $U_t$ is a
regular conditional distribution of $X$, i.e. a version of $\mathbb
P(X\in B|{\cal F}_t)=\mathbb E[I_{X^{-1}(B)}|{\cal F}_t]$. For if
$A\in{\cal F}_t$, then
\[\aligned\phantom{=}&\int_A\int_{X^{-1}(B)}\tilde{\varphi}(\omega,X(\omega'))\;R(d\omega')\mathbb
P(d\omega)\\&= \int I_{A\times X^{-1}(B)}\varphi \;d\bar{\mathbb Q}\\
&=\bar{\mathbb P}(A\times X^{-1}(B))\\&=\mathbb P(A\cap
X^{-1}(B))\\&= \int_AI_{X^{-1}(B)}\;d\mathbb P
\endaligned\]
Hence $U_t(\omega, ds) = Q_t(\omega,ds)$ $\mathbb P$--a.s. Now if
$\eta(B) = 0$, then
\[U_t(\omega,B) = \int
I_B(X(\omega')\tilde{\varphi}(\omega,X(\omega'))\;R(d\omega)=
\int_B\tilde{\varphi}(\omega,s)\;\eta(ds)=0\] so that
$Q_t(\omega,ds)<<\eta(ds)$. \eproof

\vskip0.5cm\noindent{\bf Proof of Theorem
\ref{thm_Jacod_criterion}:} Combine Theorem \ref{thm_A_implies H'}
and Proposition \ref{propn_A_is_J'}.\eproof
\begin{remarks}\rm \begin{enumerate}[(a)]\item Alternate proofs of Jacod's criterion for
$(\mathbf H')$ may be found in Jacod\cite{Jacod1985} and
Protter\cite{Protter2004}. \item Though much is made of Jacod's
criterion in the literature, it fails already in the simplest case,
when we enlarge the natural filtration of Brownian motion $W$ by
$W_T$. Indeed, $(\mathbf H')$ fails in this case. Cf. Remarks
\ref{rem_(J')_fails}.\end{enumerate}\endbox
\end{remarks}

We work in a filtered probability space $(\Omega,{\cal F},\mathbb
P,\mathbb F)$ satisfying the usual conditions.
\begin{corollary}\label{corr_enlarge_independent} Let $X$ be a random element (with state space a Borel space $(S,{\cal S})$), such that
$X$ is independent of $\mathbb F$, and let $\mathbb G$ be the
enlargement of $\mathbb F$ by $X$. Then $(\mathbf J')$ holds.
\newline Hence every $\mathbb F$--semimartingale is a $\mathbb
G$--semimartingale.
\end{corollary}
\bproof By independence, the regular conditional distributions
$Q_t(\omega, ds)$ of $X$ are equal to the law $\eta(ds)$ of $X$.
\eproof

\begin{corollary}\label{corr_enlarge_discrete} Let $X$ be a discrete random variable (i.e. takes
only countably many values), and let $\mathbb G$ be the enlargement
of $\mathbb F$ by $X$. Then $(\mathbf A)$ holds.
\newline Hence every $\mathbb F$--semimartingale is a $\mathbb
G$--semimartingale.
\end{corollary}
\bproof If $X$ has range $\{x_n:n\in\mathbb N\}$, then $\bar{\cal
F}={\cal F}_\infty\otimes\sigma(X)$ consists of countable unions of
rectangles of the form $F\times\{X=x_n\}$, where $F\in{\cal
F}_\infty$. With $\bar{\mathbb Q}=\bar{\mathbb P}|_{{\cal
F}_\infty}\otimes \bar{\mathbb P}|_{\sigma(X)}$, we see that
$\bar{\mathbb Q}(F\times\{X=x_n\}) = 0 \Rightarrow\mathbb
P(F)\mathbb P(X=x_n) = 0 \Rightarrow\mathbb P(F\cap\{X=x_n\})
=\bar{\mathbb P}(F\times\{X=x_n\}) = 0$.\eproof

\begin{corollary} {\rm (Jacod's Countable Enlargement)} Let ${\cal
A}=\{A_n:n\in\mathbb N\}$ be a family of mutually disjoint events in
${\cal F}$, and let $\mathbb G$ be the enlargement of $\mathbb F$ by
$\sigma({\cal A})$. \newline Then every $\mathbb F$--semimartingale
is a $\mathbb G$--semimartingale.
\end{corollary}
\bproof Apply Corollary \ref{corr_enlarge_discrete} to  $X = \sum_n
nI_{A_n}$.\eproof

\section{Stochastic Integrals under Enlargements}
This short section is devoted to a result of Jeulin\cite{Jeulin1980}
which gives necessary and sufficient conditions for the stochastic
integral $H\bullet M$ of a $(\mathbb F,\mathbb P)$--local martingale
w.r.t a $\mathbb F$--predictable process $H$ to be a $\mathbb
G$--semimartingale. We rely on the exposition in
Protter\cite{Protter2004}.

Recall that a process $X$ is said to be locally integrable w.r.t
$(\mathbb F,\mathbb P)$ iff there exists an increasing sequence of
$\mathbb F$--stopping times $T_n\uparrow\infty$ a.s. such that
$\mathbb E[|X_{T_n}| ;T_n>0]<\infty$ for all $n\in\mathbb N$. Recall
further that an $(\mathbb F,\mathbb P)$--semimartingale $X$ is
special iff the process $X^*_t = \sup_{s\leq t}|X_s|$ is locally
integrable. In particular, since any local martingale is obviously
special, we have that $M^*$ is locally integrable for any local
martingale $M$. If $\mathbb G$ is an enlargement of $\mathbb F$,
then every $\mathbb F$--stopping time is a $\mathbb G$--stopping
time, and hence any $(\mathbb F,\mathbb P)$--locally integrable
process is $(\mathbb G,\mathbb P)$--locally integrable. it therefore
follows that:

\begin{proposition} If an $(\mathbb F,\mathbb
P)$--local martingale is a $(\mathbb G,\mathbb P)$--semimartingale,
where $\mathbb G$ is an enlargement of $\mathbb F$, then it is a
special $(\mathbb G,\mathbb P)$--semimartingale.\endbox
\end{proposition}
\begin{theorem}\label{thm_stoch_integral_enlargement} Let $M$ be an $(\mathbb F,\mathbb P)$--local
martingale, and let $H$ be $\mathbb F$--predictable such that
$(\int_0^tH_s^2\;d[M]_s)_{t\geq 0}$ is locally integrable. Suppose
that $\mathbb G$ is an enlargement of $\mathbb F$ such that $M$
remains a $(\mathbb G,\mathbb P)$--semimartingale. Then $M$ is a
special $(\mathbb G,\mathbb P)$--semimartingale. If
$M=\tilde{M}+A$ is its $(\mathbb G,\mathbb P)$--canonical
decomposition, then $H\bullet M$ is a $(\mathbb G,\mathbb
P)$--semimartingale iff $(\int_0^tH_s\;dA_s)_{t\geq 0}$ exists as
a path--by--path Lebesgue--Stieltjes integral.\newline In that
case, the $(\mathbb G,\mathbb P)$--canonical decomposition of
$H\bullet M$ is $H\bullet M= H\bullet \tilde{M} +H\bullet A$.
\end{theorem}

\bproof LATER\dots Cf. Protter\cite{Protter2004} or
Jeulin\cite{Jeulin1980}\eproof

\section{Initial Enlargements in the Brownian World}
We follow Yor\cite{Yor1997}, Mansuy and Yor\cite{Mansuy_Yor2006}.
Let $(\Omega,{\cal F},\mathbb P)$ be a probability space that
supports a Brownian motion $W$, with natural filtration $\mathbb
F$. Let $X$ be an ${\cal F}_\infty$--measurable random variable,
and let $\mathbb G$ be the enlargement of $\mathbb F$ by $X$.
Given a bounded Borel function $f:\mathbb R\to\mathbb R$, let
$\lambda(f) = \lambda_t(f))_{t\geq 0}$ be a continuous version of
the martingale $(\mathbb E[f(X)|{\cal F}_t])_{t\geq 0}$. By the
martingale representation theorem,  there exists a unique
predictable process $\dot{\lambda}(f)$ such that
\[\lambda_t(f) = \mathbb E[f(X)]
+\int_0^t\dot{\lambda}_u(f)\;dW_u\] Denote by $\lambda_t(dx)$ the
regular conditional distributions of $X$ w.r.t. ${\cal F}_t$, i.e.
$\lambda_t(A)$ is a version of $\mathbb P(X\in A|{\cal F}_t)$. Then
\[\lambda_t(f) = \int f(x)\;\lambda_t(dx)\]
We now make the following assumption:
\begin{center}\fbox{\begin{minipage}{12cm} {\bf Assumption:} There is a family
$\dot{\lambda}_t(dx)$ of measures such that
\[\dot{\lambda}_t(f) = \int f(x)\;\dot{\lambda}_t(dx)\qquad
t\text{--a.e.}\]
\end{minipage}}\end{center}
\begin{theorem} \label{thm_Brownian_martingale_enlarged_decomp} Assume that $\dot{\lambda}_t(dx)<<\lambda_t(dx)$
$dt\times d\mathbb P$--a.e., and define $\rho(x,s)$ by
$\dot{\lambda}_t(dx) = \rho(x,t)\;\lambda_t(dx)$ Then for any
$(\mathbb F,\mathbb P)$--martingale $M =\int_0^\cdot m_u\;dW_u$
there exists a $(\mathbb G,\mathbb P)$--local martingale $\tilde
M$ such that
\[M = \tilde M +\int_0^\cdot \rho(X,u)\;d[M,W]_u\] provided that
\[\int_0^t|\rho(X,u)|\;|d[M,W]_u|<\infty \quad\text{a.s. for }t\geq
0\]
\end{theorem}
\bproof By the martingale representation theorem there is a
predictable $m$ such that $M =\int_0^\cdot m_u\;dW_u$. Let $f$ be a
bounded Borel function, and let $s<t$ and $F\in{\cal F}_s$. Then
\[\aligned &\phantom{aa}\mathbb E[I_Ff(X)(M_t-M_s)]\\
&= \mathbb E[I_F(\lambda_t(f)M_t - \lambda_s(f)M_s)]\\&=\mathbb
E\left[\int_s^tm_u\dot{\lambda}_u(f)\;du\right]\\
&=\mathbb E\left[I_F\int_s^t m_u
\left(\int\rho(x,u)f(x)\;\lambda_u(dx)\right)\;du\right]\\
&=\mathbb E\left[I_F\int_s^t m_uf(X)\rho(X,u)\;du\right]\\
&=\mathbb E\left[I_F f(X)\int_s^t\rho(X,u)\;d[M,W]_u\right]
\endaligned\]
because $d[M,W]_ = m_t\;dt$.  Hence $\mathbb E[I_Ff(X)(M_t-M_s)]
=\mathbb E\left[I_F f(X)\int_s^t\rho(X,u)\;d[M,W]_u\right]$ for all
$F\in{\cal F}_s$ and every bounded Borel function $f$. By a monotone
class argument, we see that
\[\mathbb E\left[M_t-M_s-\int_s^t\rho(X,u)\;d[M,W]_u\;\;\Big|{\cal G}_s\right] = 0\]
\eproof

Applying this result to $M=W$, and using L\'evy's characterization
of Brownian motion, we get immediately.
\begin{corollary} If, in Theorem
\ref{thm_Brownian_martingale_enlarged_decomp}, the $\mathbb
F$--Brownian motion $W$ decomposes as
\[W =\tilde W+\int_0^\cdot \rho(X,s)\;ds\]
where $\tilde W$ is a $\mathbb G$--Brownian motion (provided that
$\int_0^t|\rho(X,s)|\;ds<\infty$ a.s. for all $t\geq 0$).
\end{corollary}

The following result will allow for explicit computations:
\begin{corollary} \label{cor_enlarged_deterministic integral}Assume that $\lambda_t(dx) =\phi(t,x)\;dx$,
where $\phi(t,x)$ has the form
\[\phi(t,x) =
\phi(0,x)\exp\left(\int_0^t\rho(x,s)\;dW_s-\frac12\int_0^t\rho(x,s)^2\;ds\right)\]
Then $\dot{\lambda}_t(dx)=\rho(x,t)\;\lambda_t(dx)$, so that we
may apply Theorem \ref{thm_Brownian_martingale_enlarged_decomp}.
\end{corollary}
\bproof Since $\lambda_t(f) = \int f(x)\phi(t,x)\;dx$, we see that
$d\lambda_t(f) = \int f(x)(\phi(t,x)\rho(x,t)\;dW_t)\;dx$. Hence
\[\lambda_t(f) = \lambda_0(f) +\int_0^t\int
f(x)\phi(u,x)\rho(x,u)\;dx\;dW_u\] so that $\dot{\lambda}_u(f) =
\int f(x)\phi(u,x)\rho(x,u)\;dx$, and thus $\dot{\lambda}_t(dx) =
\phi(t,x)\rho(x,t)\;dx$. \eproof

\begin{example}\rm (Yor\cite{Yor1997}). Let $X =
\int_0^\infty\varphi(t)\;dW_t$, for some deterministic
square--integrable $\varphi$. Note that, conditional on ${\cal
F}_t$, the random variable $X$ is Gaussian with mean $m_t =
\int_0^\varphi(s)\;dW_s$ and variance $\sigma^2_t=\int_t^\infty
\varphi(s)^2\;ds$. To use Corollary
\ref{cor_enlarged_deterministic integral}, we must find
$\rho(x,s)$ so that
\[\frac{1}{\sqrt{2\pi\sigma^2_t}}\;e^{-(x-m_t)^2/2\sigma^2_t}
=\frac{1}{\sqrt{2\pi\sigma^2_0}}\;e^{-x^2/2\sigma^2_0}\;e^{\int_0^t\rho(x,s)\;dW_s-\frac12\int_0^t\rho(x,s)^2\;ds}\]
Put $M_t =\frac{1}{\sqrt{2\pi
\sigma_t^2}}e^{-(x-m_t)^2/2\sigma^2_t}$. Some messy but
straightforward calculations show that $dM_t =
M_t\frac{x-m_t}{\sigma^2_t}\varphi(t)\;dW_t$, and thus that $M_t
=e^{\int_0^t\rho(x,s)\;dW_s-\frac12\int_0^t\rho(x,s)^2\;ds}$ for
\[\rho(x,t) =\frac{x-m_t}{\sigma^2_t}\varphi(t)\]
Under suitable integrability conditions (see the Remarks that
follow), therefore, we see that $W$ is a semimartingale in the
enlargement $\mathbb G$ of $\mathbb F$ by $X$, with decomposition
\[W_t =
\tilde{W}_t+\int_0^t\frac{\varphi(s)}{\sigma_s^2}\left(\int_s^\infty\varphi(u)\;dW_u\right)\;ds\]
where $\tilde{W}$ is a $\mathbb G$--Brownian motion.\newline In
particular, if $\varphi(t) = I_{[0,T]}$, then $X =
\int_0^\infty\varphi(t)\;dW_t = W_T$, and hence
\[W_t = \tilde{W}_t +\int_0^{t\land T}\frac{W_T-W_s}{T-s}\;ds\] as
before.
\endbox
\end{example}
\begin{remarks}\label{rem_(J')_fails}\rm
Note that if $\mathbb G$ is the enlargement of $\mathbb F$ by $W_T$,
then Jacod's criterion $(\mathbf J')$ fails to hold on $[0,\infty)$.
It does, however, hold on $[0,T)$. Indeed, $W_T$ is clearly Gaussian
conditional on ${\cal F}_t$ for $t<T$. Thus if $Q_t(\omega,dx)$ are
regular conditional versions of $\mathbb P(W_T\in dx|{\cal F}_t)$,
then $Q_t(\omega,dx)$ is absolutely continuous w.r.t Lebesgue
measure. However, for $t\geq T$, $Q_t(\omega,dx)$ is the point mass
$\delta_{W_T(\omega)}(dx)$, and it is impossible to find a single
measure $\eta$ such that $\delta_{r}(dx)<<\eta(dx)$ for all
$r\in\mathbb R$.

Jeulin and Yor\cite{Jeulin_Yor1979} show that, in the Brownian
framework, a $(\mathbb F,\mathbb P)$--local martingale $M$ is a
$(\mathbb G,\mathbb P)$--semimartingale (where $\mathbb G$ is the
enlargement of $\mathbb F$ by $W_T$) iff
$\int_0^T(1-s)^{_\frac12}\;|d[M,W]_s|<\infty$, in which case
\[M_t-\int_0^{t\land T}\frac{W_1-W_s}{T-s}\;d[M,W]_s\] is a $\mathbb
G$--local martingale. Thus if $M=\int_0^\cdot m_s\;dW_s$, then $M$
remains a $\mathbb G$--semimartingale iff
$\int_0^T\frac{|m_s|}{\sqrt{T-s}}\;ds<\infty$.

Jeulin and Yor\cite{Jeulin_Yor1979} further show that there exists a
deterministic square--integrable $m_s$ such that
$\int_0^T\frac{|m_s|}{\sqrt{T-s}}\;ds=\infty$. It follows that not
every $(\mathbb F,\mathbb P)$--local martingale is a $(\mathbb
G,\mathbb P)$--semimartingale, i.e. that $(\mathbf H')$ fails.
\endbox
\end{remarks}

We follow Protter\cite{Protter2004} to elaborate on the
 results of Jeulin and Yor mentioned above:

\begin{lemma}\label{lemma_Jeulin} {\rm (Jeulin's Lemma)} On a
general filtered probability space $(\Omega,{\cal F},\mathbb
P,\mathbb F)$ satisfying the usual conditions, let
$R:\Omega\times\mathbb R^+\to\mathbb R$ be a positive measurable
process, with the properties that
\begin{enumerate}[(i)]\item all $R_s$ are identically distributed,
with common law $\mu$ satisfying $\mu\{0\}=0$ and $\int_0^\infty
x\;\mu(dx) <\infty$; \item each $R_s$ is independent of ${\cal
F}_s$.
\end{enumerate}Suppose that $A$ is a positive $\mathbb F$--predictable process with $\int_0^tA_s\;ds <\infty$ a.s. for each $t\geq 0$.
Then the following sets are $\mathbb P$--a.s.
equal:\[\left\{\int_0^\infty R_sA_s\;ds<\infty\right\}
=\left\{\int_0^\infty A_s\;ds\right\}\]\end{lemma}

\bproof (\cite{Protter2004},  Jeulin\cite{Jeulin1982}) We first show
that $\left\{\int_0^\infty R_sA_s\;ds<\infty\right\}
\subseteq\left\{\int_0^\infty A_s\;ds\right\}$. Let $E\in{\cal F}$
with $\mathbb P(E)>0$, and let $J_t = \mathbb E[I_E|{\cal F}_t]$ be
the c\`adl\`ag version of this martingale.

Let $J_*=\inf_tJ_t$. Note that $\{J_*>0\}\supseteq E$: For if
$B\subseteq E$ is an event with $\mathbb P(B)>0$, then \[\mathbb
E[J_tI_B] = \mathbb E[I_E\mathbb E[I_B|{\cal F}_t]]\geq\mathbb
E[I_B\mathbb E[I_B|{\cal F}_t]] =\mathbb E[(\mathbb E[I_B|{\cal
F}_t]^2]\geq (\mathbb E[I_B])^2 = \mathbb P(B)\] where we used the
fact that $\mathbb E\big[X\mathbb E[Y|{\cal F}_t]\big] = \mathbb
E\big[\mathbb E[X|{\cal F}_t]\;\mathbb E[Y|{\cal F}_t]\big]=\mathbb
E\big[\mathbb E[X|{\cal F}_t]Y\big]$ and also Jensen's inequality.
Now let $\{q_n:n\in\mathbb N\}$ be an enumeration of the rationals.
Then $\mathbb E[(\inf_{n\leq m} J_{q_n})I_B]\geq\mathbb P(B)$ as
well. It follows that $\mathbb E[J_*I_B]\geq \mathbb P(B)$ for all
$B\subseteq E$, and thus that $J_*>0$ a.s. on $E$.

 Let $F$ be the
common distribution function of the $R_t$, i.e. $F(x) = \mu(0,x] =
\mathbb P(R_t\leq x)$.

Note that
\[\mathbb E[I_ER_t|{\cal F}_t] =
\mathbb E\left[I_E\int_0^\infty I_{\{R_t>u\}}\;du|{\cal
F}_t\right]=\int_0^\infty\mathbb E[I_EI_{\{R_t>u\}}|{\cal
F}_t]\;du\] As $I_EI_{\{R_t>u\}} = (I_E-I_{\{R_t\leq u\}})^+$, we
have $\mathbb E[I_EI_{\{R_t>u\}}|{\cal F}_t]\geq \mathbb
E[I_E-I_{\{R_t\leq u\}}|{\cal F}_t]^+$, by Jensen's inequality.
Using the fact that $R_t$ is independent of ${\cal F}_t$, it now
follows that
\[\mathbb E[I_ER_t|{\cal F}_t]\geq
\int_0^\infty(J_t-F(u))^+\;du =\Phi(J_t)\] where \[\Phi(x):=
\int_0^\infty (x-F(u))^+\;du\] The function $\Phi$ is increasing and
continuous on $[0,1]$. Furthermore, if $x>0$, then $\Phi(x) >0$
also, as $F(0) = \mu\{0\} = 0$ by assumption.

Now replace the arbitrary event $E$ by the events \[E_n =
\left\{\int_0^\infty R_tA_t\;dt\leq n\right\}\] with $J^n_t =
\mathbb E[I_{E_n}|{\cal F}_t]$ (c\`adl\`ag), to deduce that
\[\mathbb E\left[\int_0^\infty\Phi(J_t^n)A_t\;dt\right]\leq\mathbb
E\left[\int_0^\infty\mathbb E[I_{E_n}R_t|{\cal
F}_t]A_t\;dt\right]=\mathbb E\left[I_{E_n}\int_0^\infty
R_tA_t\;dt\right]\leq n\mathbb P(E_n)<\infty\] so that
$\int_0^\infty \Phi(J_t^n)A_t\;dt<\infty$ a.s. Since
\[\Phi(J_*^n)\int_0^\infty A_t\;dt\leq\int_0\Phi(J_t^n)A_t\;dt
\qquad\text{a.s.}\] we see that \[\int_0^\infty A_t\;dt
<\infty\qquad\text{a.s. on }\;E_n\] Thus $\int_0^\infty
A_t\;dt<\infty$ a.s. on $\bigcup_nE_n = \{\int_0^\infty
R_tA_t\;dt<\infty\}$, i.e. $\{\int_0^\infty
A_t\;dt<\infty\}\supseteq\{\int_0^\infty R_tA_t\;dt<\infty\}$ a.s.

We now show the reverse inclusion: Note that if $T$ is a stopping
time, then \[\aligned \mathbb E\left[\int_0^{T}R_sA_s\;ds\right]&=
\int_0^\infty\mathbb E[I_{[0,T]}R_sA_s]\;ds\\
&=\int_0^\infty\mathbb E[I_{[0,T]}A_s\mathbb E[R_s|{\cal F}_{T\land
s}]\;ds \\&=\mathbb E\left[\int_0^TA_s\mathbb E\big[\mathbb
E[R_s|{\cal F}_s]|{\cal F}_{T\land s}\big]\;ds\right]\\
&=\alpha\mathbb E\left[\int_0^TA_s\;ds\right]
\endaligned\]
where $\alpha = \mathbb E[R_s]$ is the common mean of the $R_s$,
which is finite by assumption.

 Define the stopping time $T_n
=\inf\{t>0:\int_0^tA_s\;ds>n\}$, so that $\mathbb
E[\int_0^{T_n}A_s\;ds]\leq n$. Now
\[\mathbb E\left[\int_0^{T_n}R_sA_s\;ds\right]=
\alpha\mathbb E\left[\int_0^{T_n}A_s\;ds\right]\leq\alpha
n\<\infty\] so that $\int_0^{T_n} R_sA_s\;ds<\infty$ a.s. If
$\omega\in \{\int_0^\infty A_s\;ds<\infty\}$, then there exists $n$
(depending on $\omega$) such that $T_n(\omega) = \infty$, and so
$\int_0^\infty R_sA_s\;ds (\omega) =\int_0^{T_n}R_sA_s\;ds
(\omega)<\infty$, i.e. we also have $\omega\in\{\int_0^\infty
R_sA_s\;ds <\infty\}$. \eproof

Equipped with the lemma, we can prove the result of
Jeulin--Yor\cite{Jeulin_Yor1979}:

\begin{theorem} Let $(\Omega,{\cal F},\mathbb P)$ support a Brownian
motion $W$ with natural filtration $\mathbb F$, and let $\mathbb G$
be the enlargement of $\mathbb F$ by $W_T$. Then an $(\mathbb
F,\mathbb P)$--local martingale is an $(\mathbb G,\mathbb
P)$--semimartingale iff $\int_0^T
(T-s)^{-\frac12}\;|d[M,W]|_s<\infty$. a.s.\newline In that case
\[M_t-\int_0^{t\land T}\frac{W_T-W_s}{T-s}\;d[M,W]_s\]
is a $(\mathbb G,\mathbb P)$--local martingale.
\end{theorem}
\bproof (Protter\cite{Protter2004}) Write $M_t = M_0
+\int_0^tm_s\;dW_s$ by the martingale representation theorem,
where $m$ is predictable with $\int_0^tm_s^2\;ds<\infty$ a.s. We
know that $W$ is a $(\mathbb G,\mathbb P)$--semimartingale with
canonical decomposition $W_t =\tilde W_t -\int_0^{t\land
T}\frac{W_T-W_s}{T-s}\;ds$. Using Theorem
\ref{thm_stoch_integral_enlargement}, we see that $M$ is a
$(\mathbb G,\mathbb P)$--semimartingale iff
$\int_0^t|m_s|\frac{|W_T-W_s|}{T-s}\;ds<\infty$ a.s. for $0\leq
t\leq T$. Now apply Jeulin's lemma, with $A_t
=\frac{|m_t|}{\sqrt{T-t}}$ and $R_t =
I_{\{t<T\}}\frac{|W_T-W_t|}{\sqrt{T-t}}$. (Note that for $t<T$, we
have $R_t$ independent of ${\cal F}_t$, with law $N(0,1)$, so that
Jeulin's Lemma is applicable.) Then
$\int_0^t|m_s|\frac{|W_T-W_s|}{T-s}\;ds = \int_0^tA_sR_s\;ds$
which is finite a.s. iff
$\int_0^tA_s\;ds=\int_0^t\frac{m_s}{\sqrt{1-s}}\;ds$ is finite
a.s. Now by the associative law for integrals,
\[\int_0^t\frac{|m_s|}{\sqrt{1-s}}\;ds
=\int_0^t\frac{1}{\sqrt{T-s}}\;|d[m\bullet W,W]_s|
=\int_0^t\frac{1}{\sqrt{T-s}}\;|d[M,W]_s|\]

Again by Theorem  \ref{thm_stoch_integral_enlargement}, we see
that $(m\bullet W)_t - \int_0^{T\land t}
\frac{m_s(W_T-W_s)}{T-s}\;ds = M_t-\int_0^{T\land
t}\frac{W_T-W_s}{T-s}\;d[M,W]_s$ is a $(\mathbb G,\mathbb
P)$--local martingale.\eproof
\begin{example}\rm Here is an example of Jeulin and
Yor\cite{Jeulin_Yor1979} of a $(\mathbb F,\mathbb P)$--local
martingale which is not a $(\mathbb G,\mathbb P)$--semimartingale
(where $\mathbb F$ is the natural filtration of Brownian motion
$W$, and $\mathbb G$ is the enlargement of $\mathbb F$ with
$W_1$): Take an $\alpha\in(\frac12,1)$, and let $m_s
=(1-s)^{-\frac12}(-\ln(1-s))^{-\alpha}I_{\{\frac12<s<1\}}$. Then
$m$ is a deterministic predictable function with
$\int_0^1m_s^2\;ds <\infty$, so that $M=m\bullet W$ is defined.
However, $\int_0^1\frac{m_s}{\sqrt{1-s}} =\infty$, and so $M$ is a
$(\mathbb F,\mathbb P)$--local martingale, but not a $(\mathbb
G,\mathbb P)$--semimartingale.\endbox
\end{example}

\newpage
\section{References}
\bibliographystyle{plain}
\bibliography{Bib_Enlargement_Notes}

\begin{thebibliography}{10}

\bibitem{Ankirchner2005}
Stefan Ankirchner.
\newblock {\em Information and Semimartingales}.
\newblock PhD thesis, Humboldt University, Berlin, 2005.

\bibitem{ADI2005}
Stefan Ankirchner, Steffen Dereich, and Peter Imkeller.
\newblock Enlargement of filtrations and continuous girsanov-type embeddings.
\newblock {\em preprint: http://www.mathematik.hu-berlin.de/~ankirchn/}, 2005.

\bibitem{Delbaen_Schachermayer1995}
Freddy Delbaen and Walter Schachermayer.
\newblock A general version of the fundamental theorem of asset pricing.
\newblock {\em Mathematische Annalen}, 300:463--520, 1995.

\bibitem{Dellacherie_Meyer1978}
Claude Dellacherie and Paul-Andr\'e Meyer.
\newblock {\em Probabilities and Potential A}.
\newblock North--Holland, 1978.

\bibitem{Jacod1985}
Jean Jacod.
\newblock Grossissement initial, hypoth\'ese h', et th\'eor\`eme de girsanov.
\newblock In Thierry Jeulin and Marc Yor, editors, {\em Grossisements de
  filtrations: exemples et applications}, LNM 1118, pages 15--35. Springer,
  1985.

\bibitem{Jeulin1980}
Thierry Jeulin.
\newblock {\em Semi-Martingales et Grossissement d'une Filtration}.
\newblock LNM 833. Springer, 1980.

\bibitem{Jeulin1982}
Thierry Jeulin.
\newblock Sur la convergence absolue de certaines int\'egrales.
\newblock In {\em S\'eminaire de probabilit\'es (Strasbourg) XVI}, LNM 920,
  pages 248--256. Springer, 1982.

\bibitem{Jeulin_Yor1979}
Thierry Jeulin and Marc Yor.
\newblock In\`egalit\`e de hardy, semimartingales, et faux amis.
\newblock In {\em S\'eminaire de probabilit\'es (Strasbourg) XIII}, LNM 721,
  pages 332--359. Springer, 1979.

\bibitem{Kallenberg2002}
Olav Kallenberg.
\newblock {\em Foundations of Modern probability}.
\newblock Springer, 2nd edition, 2002.

\bibitem{Karatzas_Shreve1991}
Ioannis Karatzas and Steven~E. Shreve.
\newblock {\em Brownian Motion and Stochastic Calculus}.
\newblock Springer, 2$^{\text{nd}}$ edition, 1991.

\bibitem{Mansuy_Yor2006}
Roger Mansuy and Marc Yor.
\newblock {\em Random Times and Enlargements of Filtrations in a Brownian
  Setting}.
\newblock LNM 1873. Springer, 2006.

\bibitem{Parthasarathy1967}
K.R. Parthasarathy.
\newblock {\em Probability measures on metric spaces}.
\newblock Academic Press, 1967.

\bibitem{Protter2004}
Philip Protter.
\newblock {\em Stochastic Integration and Differential Equations}.
\newblock Springer, 2004.

\bibitem{Yor1997}
Marc Yor.
\newblock {\em Some Aspects of Brownian Motion, Part II: Some Recent Martingale
  Problems}.
\newblock Lectures in Mathematics ETH. Birkh\"auser, 1997.

\bibitem{YM1978}
Marc Yor and Paul-Andr\'e Meyer.
\newblock Sur l'extension d'un th\'eor\`eme de doob \`a un noyau
  $\sigma$--fini, d'apr\`es g.mokobodzki.
\newblock In {\em S\'eminaire de Probabilit\'es (Strasbourg) XII}, LNM 649,
  pages 482--488. Springer, 1978.

\end{thebibliography}

\newpage

\appendix
\section{Monotone Class Theorems}
The monotone class theorems are a collection of related results that
prove that if a certain \lq\lq nice" set of measurable functions
satisfy a certain property, then all (bounded) measurable function
have that property. In the literature, one frequently finds that a
proof will verify a property for indicator functions, and then
assert that \lq\lq the rest follows by a monotone class
argument".

The monotone class theorems are based on a kind of \lq\lq
decomposition" of a $\sigma$--algebra into a part that meshes nicely
with the properties of measures ($\lambda$--systems) and part which
doesn't, but which is nevertheless very simple, and meshes nicely
with multiplication ($\pi$--systems):
\begin{definition}\label{defn_pi_lambda}\rm Let ${\cal C}$ be a collection of subsets of $\Omega$
\begin{enumerate}[(a)] \item ${\cal C}$ is called a
$\pi$--system if it is closed under finite intersections. \item
${\cal C}$ is called a $\lambda$--system if
\begin{enumerate}[(i)]\item $\Omega\in{\cal C}$;
\item $A, B\in {\cal C}$ and $A\subseteq B$ implies $B-A\in{\cal
C}$; \item If $A_1,A_2,\dots\in{\cal C}$ and $A_n\uparrow A$, then
$A\in{\cal C}$.
\end{enumerate}
\end{enumerate}
\endbox
\end{definition}
We denote by $\pi({\cal C})$ and $\lambda({\cal C})$ the $\pi$--,
respectively, $\lambda$--system {\em generated} by ${\cal C}$ The
following lemma is easy to prove.
\begin{lemma}\label{lemma_pi_lambda} A family ${\cal C}$ of subsets of $\Omega$
is a $\sigma$--algebra iff it is {\em both} a $\pi$--system and a
$\lambda$--system.\endbox
\end{lemma}

The following technical result often allows us to work with \lq\lq
easy" $\pi$--systems, instead of the \lq\lq difficult"
$\sigma$--algebras:
\begin{theorem}\label{thm_pi_lambda} {\rm (Dynkin's
Lemma, Monotone Class Theorem)}\begin{enumerate}[(a)] \item  If
${\cal C}$ is a $\pi$--system on $\Omega$, then \[\lambda({\cal
C})=\sigma({\cal C})\] \item Suppose that ${\cal C}$ is a
$\pi$--system
 and that ${\cal D}$ is a $\lambda$--system
(both on a set $\Omega$), and also that ${\cal C}\subseteq{\cal D}$.
Then $\sigma({\cal C})\subseteq{\cal D}$.
\end{enumerate}
\end{theorem}
\bproof (a) Let ${\cal D} = \lambda({\cal C})$. By Lemma
\ref{lemma_pi_lambda}, it suffices to show that ${\cal D}$ is a
$\pi$--system. We do this in two steps.\\STEP I: Fix $C\in{\cal C}$,
and define \[{\cal D}_C =\{A\in{\cal D}: A\cap C\in{\cal D}\}\] Then
${\cal C}\subseteq{\cal D}_C\subseteq {\cal D}$ (because ${\cal C}$
is a $\pi$--system). Then ${\cal D}_C$ is easily shown to be a
$\lambda$--system containing ${\cal C}$, so that ${\cal D}_C = {\cal
D}$. \newline STEP II: Now, fix any $D\in{\cal D}$, and define
\[{\cal D}^D=\{A\in{\cal D}: A\cap D\in{\cal D}\}\]
First note that if $C\in{\cal C}$, then ${\cal D}_C = {\cal D}$, so
$D\in{\cal D}_C$. It follows that $D\cap C\in{\cal D}$, and thus
that $C\in{\cal D}^D$, for every $C\in{\cal C}$. Thus ${\cal
C}\subseteq {\cal D}^D$, for all $D\in{\cal D}$.\\It follows as
above that ${\cal D}^D$ is a $\lambda$--system, and thus that ${\cal
D}^D={\cal D}$, for all $D\in{\cal D}$.

In particular, if $A,B\in{\cal D}$, then $A\in{\cal D}^B$, and so
$A\cap B\in{\cal D}$. This shows that ${\cal D}$ is a $\pi$--system,
and thus a $\sigma$--algebra.

(b) follows directly from (a). \eproof

\begin{definition}\rm
\begin{enumerate}[(a)]
\item A collection ${\cal A}$ of bounded real--valued functions on a
set $\Omega$ is called an {\em algebra}\footnote{This is not to be
confused with a family of sets closed under complementation and
finite unions/intersections.} if it is a vector space and closed
under multiplication.
\item A collection ${\cal H}$ of bounded real--valued functions on a
set $\Omega$ is called a {\em monotone vector space} iff
\begin{enumerate}[(i)]\item ${\cal H}$ is a vector space over
$\mathbb R$.
\item The constant function 1 belongs to ${\cal H}$.
\item If $(f_n)_n$ is a uniformly bounded increasing sequence of non--negative members of
${\cal H}$, then $\lim_nf_n\in{\cal H}$.

\end{enumerate}
\end{enumerate}
\endbox
\end{definition}

\begin{theorem} \label{thm_monotone_class}{\rm (Monotone Class Theorem)}\newline
Let ${\cal H}$ be a monotone vector space on $S$. Let ${\cal A}$ be
a $\pi$--system on $S$ with the property that $I_A\in{\cal H}$ for
every $A\in{\cal A}$.
\newline
Then every bounded $\sigma({\cal A})$--measurable function belongs
to ${\cal H}$.
\end{theorem}
\noindent{\bf Proof:}
 Let ${\cal D} = \{F\subseteq S: I_F\in{\cal H}\}$. It is not hard to
 show that ${\cal D}$ is a $\lambda$--system.
By Theorem \ref{thm_pi_lambda},  ${\cal D}\supseteq\sigma({\cal
A})$.

Let $h$ be a non--negative, bounded $\sigma({\cal A})$--measurable
function, with upper bound $K$, i.e.
\begin{equation}
0\leq h(s)\leq K\qquad\text{for all }s\in S \notag
\end{equation}
If we define
\[
h_n(s) = \sum\limits_{k=1}^{K2^n}\frac{k-1}{2^n}I_{A(n,k)}(s)
\qquad\text{where}\qquad A(n,k) = \left\{s\in S:\frac{k-1}{2^n}\leq
h(s)<\frac{k}{2^n}\right\}
\]

then the $h_n$ are simple functions with $h_n\uparrow h$. Since $h$
is $\sigma({\cal A})$--measurable, each $A(n,k)\in{\cal D}$, i.e.
$I_{A(n,k)}\in{\cal H}$. Because ${\cal H}$ is a vector space, we
now see that $h_n\in{\cal H}$ for each $n\in\mathbb N$. Thus
$h\in{\cal H}$ as well.

We have now shown that every non--negative bounded $\sigma({\cal
A})$--measurable function belongs to ${\cal H}$. The same result can
be obtained for arbitrary bounded $h$ by splitting into positive and
negative parts: $h = h^+-h^-$. \eproof

Here is another such result:

\begin{theorem}\label{thm_monotone_class_2} {\rm (Monotone Class Theorem)}  Let ${\cal M}$ be a collection of bounded
real--valued functions on $\Omega$ which is closed under
multiplication.  Suppose that ${\cal H}$ is a monotone vector space
which is closed under uniform convergence, such that ${\cal
H}\supseteq{\cal M}$. Then every bounded $\sigma({\cal
M})$--measurable function belongs to ${\cal H}$.
\end{theorem}
\bproof We sketch the proof given in Dellacherie and
Meyer\cite{Dellacherie_Meyer1978}. We may assume that $1\in{\cal
M}$. Let ${\cal A}$ be an  algebra which is maximal such that ${\cal
M}\subseteq{\cal A}\subseteq{\cal H}$ --- such ${\cal A}$ exists by
Zorn's Lemma. The function $x\mapsto|x|$ can be uniformly
approximated by polynomials on every compact interval of $\mathbb
R$. (To see this note that $|x|=(1-(1-x^2))^\frac12$, and that the
Taylor series of $z\mapsto (1-z)^\frac12$ converges uniformly on
$[-1.1]$.) Now ${\cal A}$ is closed under uniform convergence
because ${\cal H}$ is, by the maximality of ${\cal A}$. It follows
that if $f\in{\cal A}$, then also $|f|\in{\cal A}$. In particular,
given $f,g\in{\cal A}$, we see that $f^\pm =\frac{|f|\pm
f}{2}\in{\cal A}$, and hence that $f\lor g = g+(f-g)^+\in{\cal A}$
and  $f\land g = f+g-f\lor g\in{\cal A}$.

We claim that ${\cal A}$ is a monotone vector space. For suppose
that $g$ is the limit of a uniformly bounded increasing sequence
$(g_n)_n$ of members of ${\cal A}$.  It is easy to see that the
algebra generated by ${\cal A}\cup\{g\}\subseteq{\cal H}$, so that
$g\in{\cal A}$ by maximality of ${\cal A}$.

Let ${\cal C} = \{C\subseteq\Omega: I_C\in{\cal A}\}$. Since ${\cal
A}$ is an algebra and $1\in{\cal A}$, we see that ${\cal C}$ is a
closed under finite intersections and complementation, and thus also
under finite unions. Since ${\cal A}$ is monotone, ${\cal C}$ is
closed under countable unions, and thus a $\sigma$--algebra.

Since every non--negative bounded measurable function can be
uniformly approximated from below by simple measurable functions (as
in the proof of Theorem \ref{thm_monotone_class}), we see that
${\cal A}$ contains all ${\cal C}$--measurable functions.

To complete the proof, we need only show that $\sigma({\cal
M})\subseteq{\cal C}$, i.e. that $\{f\leq c\}\in{\cal C}$ for all
$f\in{\cal M}$ and all $c\in\mathbb R$. Clearly it suffices to show
that $\{f\geq 1\}\in{\cal C}$ for all $f\in{\cal A}$, i.e. that
$I_{\{f\geq 1\}}\in{\cal A}$ for all $f\in{\cal A}$. Now if
$f\in{\cal A}$, then $g=(f\land 1)^+\in{\cal A}$. Now the sequence
of $n^{\text{th}}$ powers of $g$ has $g^n\downarrow I_{\{f\geq
1\}}$, so that $I_{\{f\geq 1\}}\in{\cal A}$ because ${\cal A}$ is a
monotone algebra. \eproof

\begin{remarks}\rm Protter(\cite{Protter2004}, p.7) states that we may drop
the assumption that ${\cal H}$ is closed under uniform convergence
in Theorem \ref{thm_monotone_class_2}: A monotone vector space is
always closed under uniform convergence.
\endbox
\end{remarks}

\section{Completions, the Usual Hypotheses, etc.} Here
follow some definitions and remarks:
\begin{enumerate}[1.]\item
Let $(\Omega,{\cal F},\mathbb P)$ be a probability space. ${\cal F}$
is said to be complete iff it contains all $\mathbb P$--negligible
sets: Whenever $A\subseteq B$ and $B\in{\cal F}$ with $\mathbb P(B)
= 0$, then $A\in{\cal F}$.
\item
Suppose that $(\Omega,{\cal F},\mathbb P)$ is a complete probability
space, and that ${\cal G}\subseteq{\cal F}$. Then the completion
${\cal G}^\mathbb P$ of ${\cal G}$ in ${\cal F}$ is the
$\sigma$--algebra with the following property: \[A\in{\cal
G}^\mathbb P\qquad\text{ iff there exists }B\in{\cal F}\text{ such
that }\quad \mathbb P(A\Delta B) = 0\] Equivalently,
\[A\in{\cal G}^{\mathbb P}\quad\text{ iff there exists }B\in{\cal
G}\text{ and null sets }M,N\in{\cal F}\text{ such that }B-N\subseteq
A\subseteq B\cup M\]
\item
A filtration $\mathbb F=({\cal F}_t)_{t\geq 0}$ is right--continuous
if and only if \[{\cal F}_t = {\cal F}_{t+}:=\bigcap_{s>t}{\cal
F}_s\qquad\text{for all }t>0\]
\item A filtration $\mathbb F=({\cal F}_t)_{t\geq 0}$ on a complete
probability space $(\Omega,{\cal F},\mathbb P)$ is said to satisfy
the {\em usual conditions} (w.r.t a probability measure $\mathbb P$)
iff
\begin{enumerate}[(i)]\item $\mathbb F$ is right-- continuous; and,
\item ${\cal F}_0$ contains all the $\mathbb P$--null sets in ${\cal
F}$.
\end{enumerate}
\item If ${\cal F}$ is a filtration, we can augment it to satisfy
the usual conditions, as follows: Let
\[{\cal G}_t =\bigcap_{s>t}({\cal F}_s\cap\sigma({\cal N}))\] where
${\cal N}$ is the collection of all $\mathbb P$--null sets.
\item The imposition of the usual conditions is essential for the existence
of regular versions of stochastic processes: Recall that a $(\mathbb
F,\mathbb P)$-- submartingale $X=(X_t)_{t\geq 0}$has a c\`adl\`ag
version iff the map $t\mapsto\mathbb EX_t$ is right--continuous.
However, this requires the usual conditions; cf. Karatzas and
Shreve\cite{Karatzas_Shreve1991}. To be precise, if $X=(X_t)_t$ is
an $(({\cal F}_t)_t,\mathbb P)$--submartingale (not assuming the
usual conditions on $({\cal F}_t)_t$), then $X^+=(X_{t+})_t$ is a
c\`adl\`ag $(({\cal F}_{t+})_t,\mathbb P)$--submartingale. Here
$X_{t+} = \lim_{s\downarrow t} X_s$ exists for all $t$, a.s. If
$({\cal F}_t)_t$ is right--continuous, then $X^+$ is adapted to
$({\cal F}_t)_t$, and can then be shown to be a modification of $X$.
Cf. \cite{Karatzas_Shreve1991} p.16 for more details.
\end{enumerate}

\section{Regular Conditional Probabilities}

Let $(\Omega,{\cal F},\mathbb P)$ be a probability space, and let
$(S,{\cal S}), (T,{\cal T})$ be measurable spaces. Given  random
elements $X:(\Omega,{\cal F},\mathbb P)\rightarrow(S,{\cal S})$
and $Y:(\Omega,{\cal F},\mathbb P)\rightarrow(T,{\cal T})$,  we
know how to define $\mathbb P(X\in A|Y)$ for $A\in {\cal S}$:\[
\mathbb P(X\in A|Y) = \mathbb E[I_A(X)|\sigma(Y)]\] is just a
version of a particular conditional expectation (which is a.s.
unique).

We want, however, to make sense of the expression
\[\mathbb P(X\in A|Y=y)\qquad \text{for }A\in{\cal S}, y\in T\]
Consider the joint law $\mathbb P_{X,Y}$ on $(S\times T,{\cal
S}\otimes{\cal T})$, given by
\[\mathbb P_{X,Y}(A\times B) =\mathbb P(X\in A,Y\in B)\]
We ought then be able to write
\[\mathbb P_{X,Y}(A\times B) = \int_B\mathbb P(X\in A|Y=y)\;\mathbb
P_Y(dy)=\int_B\mathbb P_X^y(A)\;\mathbb P_Y(dy)\] where $\mathbb
P_Y$ is the law of $Y$ on $(T,{\cal T})$, and $\mathbb P_X^y(A) =
\mathbb P(X\in A|Y=y)$. If we can do this, we have {\em
disintegrated} the joint law. However, to be able to do it, we
clearly require that
\begin{itemize}\item For $y\in T$, each map
$P_X^y$ is a probability measure on $(S,{\cal S})$, and\item For
fixed $A\in{\cal S}$, the map $y\mapsto \mathbb P(X\in A|Y=y)$ is
measurable, so that we can perform the integration.
\end{itemize}
\begin{definition}\rm Given two measurable spaces $(S,{\cal S}),
(T,{\cal T})$, a map $\mu:T\times{\cal S}\to\bar{\mathbb R}^+$ is
called a {\em stochastic kernel} from $(T,{\cal T})$ to $(S,{\cal
S})$ iff
\begin{enumerate}[(i)]\item The map $t\mapsto\mu(A,t)$ is ${\cal
T}$--measurable in $t\in T$ for fixed $A\in{\cal S}$, and \item
The map $A\mapsto\mu(A,t)$ is a probability measure on $(S,{\cal
S})$ for fixed $t\in T$.
\end{enumerate}
If $(T,{\cal T})=(\Omega,{\cal F})$ is a probability space, then
the stochastic kernel $\mu$ is called a {\em random
measure}.\endbox
\end{definition}

\begin{definition}\label{Defn_regular_probability}\rm Let $(\Omega,{\cal F},\mathbb P)$ be a probability space. Given  random
elements $X:(\Omega,{\cal F},\mathbb P)\rightarrow(S,{\cal S})$
and $Y:(\Omega,{\cal F},\mathbb P)\rightarrow(T,{\cal T})$, a {\em
regular conditional probability of $X$ given $Y$} is a random
measure of the form \[\mu(Y, A) = \mathbb P[X\in
A|Y]\quad\text{a.s.}\] where $\mu$ is a stochastic kernel from $T$
to $S$. \newline We then define $\mathbb P[X\in A|Y=y] :=
\mu(y,A)$.\endbox
\end{definition}

\begin{remarks}\rm Given $X:(\Omega,{\cal F},\mathbb P)\mapsto (S,{\cal
S})$ and a sub--$\sigma$--algebra${\cal G}$  of ${\cal F}$, one
often encounters the notion of a regular conditional probability
of the form $\mathbb P[X\in \cdot|{\cal G}]$. By this is meant a
version of $\mathbb P[X\in \cdot|{\cal G}]$ which is a stochastic
kernel from $(\Omega,{\cal F})$ to $(S,{\cal S})$, i.e. a map
$\nu:\Omega\times{\cal S}\to \bar{\mathbb R}^+$
having\begin{enumerate}[(i)]\item $A\mapsto\nu(\omega,A)$ is a.s.
a probability measure on $(S,{\cal S})$; \item
$\omega\mapsto\nu(\omega,A)$ is ${\cal F}$--measurable ---
actually, ${\cal G}$--measurable, because in addition\item
$\nu(\cdot, A)$ is a version of $\mathbb E[X\in A|{\cal
G}](\cdot)$.
\end{enumerate}
The preceding Defn \ref{Defn_regular_probability} includes this as
a special case: Put $Y = \text{Id}_\Omega:(\Omega,{\cal
F})\to(\Omega,{\cal G})$, and Then $\sigma(\text{Id}) = {\cal G}$,
so that $\mu(\omega, A)=\mu(Y,A)(\omega) = \mathbb P[X\in A|{\cal
G}]$ a.s.
\endbox
\end{remarks}
Regular conditional probabilities do not always exist --- cf. Rogers
and Williams I.43 for a counterexample --- but some mild topological
conditions on the state space $(S,{\cal S})$ will ensure existence:
A {\em Borel space} is a measurable space $(S,{\cal S})$ with the
property that there exists a Borel subset $B\subseteq[0,1]$ and a
bijection $f:(S,{\cal S})\to (B,{\cal B}(B))$ such that both $f,
f^{-1}$ are measurable. In particular, it is known that any Polish
space (equipped with its Borel algebra) is a Borel space (cf.
Parthasarathy\cite{Parthasarathy1967}). Consequently, any Borel
subset of a Polish space, is a Borel space.
\begin{theorem} \label{thm_regular_probabilities_exist}{\rm (Existence of Regular Conditional
Distributions)} Let $(\Omega,{\cal F},\mathbb P)$ be a probability
space, and let $(S,{\cal S}), (T,{\cal T})$ be measurable spaces,
and assume in addition that $(S,{\cal S})$ is a Borel space. Given
random elements $X:(\Omega,{\cal F},\mathbb P)\rightarrow(S,{\cal
S})$ and $Y:(\Omega,{\cal F},\mathbb P)\rightarrow(T,{\cal T})$.
Then there exists a unique regular conditional distribution of $X$
given $Y$ (i.e. there exists probability kernel $\mu:T\times{\cal
S}\to\bar{\mathbb R}^+$ such that $\mathbb P[X\in\cdot|Y] =
\mu(Y,\cdot)$ a.s., and any two such kernels are $\mathbb
P_Y$--a.e. equal.)
\end{theorem}

\bproof (Taken from Kallenberg\cite{Kallenberg2002}.) Since $S$ is
a Borel space, we may assume w.l.o.g. that $S\in{\cal B}(\mathbb
R)$, and thus that $X$ is real--valued, and ${\cal S} = S\cap{\cal
B}(\mathbb R)$.

 Choose, for
$q\in\mathbb Q$, a measurable function $f_q:T\to[0,1]$ such that
\[f_q(Y) = \mathbb P[X\leq q|Y]\text{ a.s.}\]
(this is possible, by the Doob--Dynkin Lemma) and define $f(t,q) =
f_q(t)$. Let \[T'=\{t\in T: f(t,q)\text{ is increasing in }q,
\lim\limits_{q\to\infty}f(t,q)=1, \lim\limits_{q\to -\infty}f(t,q)
= 0\}\] Note that $T'\in{\cal T}$ (because $T'= \bigcap_{q_1\leq
q_2}\{t:f_{q_1}(t)\leq
f_{q_2}(t)\}\cap\bigcap_n\bigcup_q\{t:f_q(t)<\frac1n\}\cap\bigcap_n\bigcup_q\{t:f_q(t)>1-\frac1n\}$).
Also note that each of the conditions defining $T'$ holds at
$Y(\omega)$ a.s., so that $Y\in T'$ a.s.

We now define $F:T\times\mathbb R\to[0,1]$ by
\[F(t,x)= I_{T'}\inf_{q>x}f(t,q) +I_{T'^c}I_{\{x\geq 0\}}\] which
has the property that $F(t,\cdot)$ is a distribution function on
$\mathbb R$ for all $t\in\mathbb T$ (recall that a function is a
distribution function precisely when it is right--continuous,
tends to 1 at $+\infty$ and to 0 at $-\infty$). Let $m_t$ be the
Lebesgue--Stieltjes measure associated with $F(t,\cdot)$ (i.e. the
unique probability measure on $(\mathbb R,{\cal B}(\mathbb R)$
satisfying $m_t(-\infty, x] = F(t,x)$). Note also that, for fixed
$x$, the map $t\mapsto F(t,x)$ is measurable in $t\in T$ (because,
e.g. for $0<u<1$ and $x\geq 0$, we have $\{t: F(t,x)<u\} = \{t\in
T': \exists q\in \mathbb Q(q>x\land f(t,q)<u)\}$, etc.) So $m(t,
B) =m_t(B)$ behaves like a stochastic kernel on the sets
$B=(-\infty,x]$, and a monotone class argument shows that $m$ is a
stochastic kernel from $(T,{\cal T})$ to $(\mathbb R,{\cal
B}(\mathbb R))$.

Furthermore, since $m(t,(-\infty,x]) = F(t,x)$, we have
$m(Y,(-\infty,x]) = F(Y,x)$. Since also $f(Y,q) = \mathbb P[X\leq
q|Y]$ a.s., the monotone convergence theorem for conditional
expectations ensures that $F(Y,x) = \mathbb P[X\leq x|Y]$ a.s., so
that $m(Y,(-\infty x]) = \mathbb P[X\leq x|Y]$ a.s. Another
application of a monotone class theorem shows that
\[m(Y,B) = \mathbb P[X\in B|Y]\text{ a.s.}\qquad\text{for all }
B\in{\cal B}(\mathbb R)\] The kernel $m$ is almost what we seek:
it is a kernel form $T$ to $\mathbb r$, whereas we need a kernel
from $T$ to $S$. Note, however, that $m(Y,S^c) = 0$ a.s. (because
$X$ takes values in $S$). Define $\mu:T\times{\cal S}\to [0,1]$ by
\[\mu(t,\cdot) =\left\{\aligned
m(t,\cdot)\qquad&\text{if }m(t,S)=1\\
\delta_{s_0}\qquad&\text{else}\endaligned\right.\] where $s_0\in
S$ is arbitrary. It is not hard to see that $\mu$ is a regular
conditional distribution.

Uniqueness is easy: If $\mu'$ is another regular conditional
distribution of $X$ given $Y$, then since two versions of
conditional expectation are a.s. equal, we have
\[\mu(Y,(-\infty,q]) = \mathbb P[X\leq q|Y] =\mu'(Y,(-\infty,q])
\text{ a.s.}\] By a monotone class theorem, $\mu(Y,\cdot) =
\mu'(Y,\cdot)$.

\eproof

\begin{theorem}\label{Thm_disintegration} {\rm (Disintegration)} Let $(\Omega,{\cal F},\mathbb P)$ be a probability
space, and let $(S,{\cal S}), (T,{\cal T})$ be measurable spaces.
Suppose we are given a sub--$\sigma$--algebra ${\cal G}$ of ${\cal
F}$ and a random element $X$ of $S$ such that $\mathbb
P[X\in\cdot|{\cal G}]$ has a regular conditional version $\nu$.
Also, consider a ${\cal G}$--measurable random element $Y$ of $T$
and a measurable function $f:S\times T\to\mathbb R$ with $\mathbb
E|f(X,Y)|<\infty$. Then
\[\mathbb E[f(X,Y)|{\cal G}] =\int_S
f(s,Y)\;\nu(ds)\quad\text{a.s.}\tag{$\dagger$}\]
\end{theorem}
\bproof (Taken from Kallenberg\cite{Kallenberg2002}.) If $f =
I_{A\times B}$, where $A\in{\cal S}, B\in{\cal T}$, then $\mathbb
P[X\in A,Y\in B] =\mathbb E[I_B(Y) \mathbb E[I_A(X)|{\cal
G}]\;]=\mathbb E[I_B(Y)\int_SI_A(s)\;\nu(ds)] =\mathbb
E\left[\int_Sf(s,Y)\;\nu(ds)\right]$. Thus \[\mathbb
E[f(X,Y)]=\mathbb E\left[\int_Sf(s,Y)\;\nu(ds)\right]\tag{$*$}\]
for $f=I_{A\times B}$. By a monotone class theorem, $(*)$ holds
for all measurable indicator functions, then extends by linearity
and the monotone convergence theorem to all non--negative
measurable functions.

Now assume that $f\geq 0$ is measurable with $\mathbb
E[f(X,Y)<\infty$, and fix $G\in{\cal G}$. Then $(Y,I_A)$ is a
${\cal G}$--measurable random element of $T\times\{0,1\}$, so
using $(*)$ we see that
\[\mathbb E[f(X,Y)I_A] = \mathbb E\left[\int_Sf(s,Y)I_A\;\nu(ds)\right] =
\mathbb E\left[\int_Sf(s,Y)\;\nu(ds)\;I_A\right]\] from which it
follows that $\int_Sf(s,Y)\;\nu(ds)$ is a version of $\mathbb
E[f(X,Y)|{\cal G}]$. This proves the result for $f\geq 0$. The
general result follows by decomposing a measurable $f$ into a
difference of its positive and negative parts.\eproof

Taking ${\cal G} = \sigma(Y)$ in $(\dagger)$ in Thm.
\ref{Thm_disintegration}, and $\nu(\omega, ds) =
\mu(Y(\omega),ds)$ for we see that
\begin{corollary} If $\mathbb P[X\in\cdot|Y]$ has a regular
conditional distribution $\mu$, then
\[\mathbb E[f(X,Y)|Y] =\int_Sf(s,Y)\;\mu(Y, ds)\]
whenever $f$ is measurable with $\mathbb E|f(X,Y)|<\infty$.\endbox
\end{corollary}
Applying the expectation operator to both sides of the equation in
the preceding corollary, we see that
\begin{corollary}If $\mathbb P[X\in\cdot|Y]$ has a regular
conditional distribution $\mu$, then
\[\mathbb E[f(X,Y)] =\mathbb E\left[\int_Sf(s,Y)\;\mu(Y,ds)\right]\]
whenever $f$ is measurable with $\mathbb E|f(X,Y)|<\infty$.\endbox
\end{corollary}

We end this section with a result on \lq\lq nice", jointly
measurable, densities. Recall that a $\sigma$--algebra is said to be
{\em separable} if it is generated by a countable family of sets
(i.e. ${\cal F}$ is separable iff ${\cal F} = \sigma(F_n:n\in\mathbb
N)$). The Borel algebra of a separable metrizable space $X$ is
clearly separable: For if $D$ is a countable dense subset of $X$,
then ${\cal B}(X) =\sigma(B(d,\frac1n): d\in D,n\in\mathbb N)$.
\begin{theorem}\label{thm_Doob_disintegration} {\rm (Doob Theorem on Disintegration)}
Let $\mu:\Omega\times{\cal S}\to\bar{\mathbb R}^+$ be a  stochastic
kernel from $(T,{\cal T})$ to $(S,{\cal S})$, and let $\eta$ be a
finite measure on $(S,{\cal S})$. Suppose in addition that:
\begin{enumerate}[(i)]\item ${\cal S}$ is separable.
\item Each measure $\mu(t,\cdot)$ is absolutely continuous w.r.t.
$\eta$ on $(S,{\cal S})$.
\item Each measure $\mu(t,\cdot)$ is finite on $(S,{\cal
S})$, for all $t\in T$.\end{enumerate} Then there is a non--negative
${\cal T}\otimes{\cal S}$--measurable function $g(t,s)$ such that
\[ \mu(t,ds) = g(t,s)\;\eta(ds)\]
\end{theorem}
\bproof (Taken from Yor and Meyer\cite{YM1978}) Without loss of
generality, we may assume that $\eta =\mathbb P$ is a probability
measure. As ${\cal S}$ is separable, we may write ${\cal
S}=\bigvee_n{\cal S}_n$, where the ${\cal S}_n$ form a sequence of
sub--$\sigma$--algebras of ${\cal S}$ generated by finer and finer
finite partitions ${\cal P}_n$. Define
\[g_n(t,s) = \frac{d\mu(t,ds)}{d\mathbb P(ds)}\Big|_{{\cal S}_n}\]
i.e. if $s\in A$, where $A\in{\cal P}_n$ is a block of the partition
that generates ${\cal S}_n$, then $g^n(t,s) =\mu(t,A)/\mathbb P(A)$.
That
 $g_n(t,s)$ is ${\cal T}\otimes{\cal
S}_n$--measurable now follows from the measurability of
$t\mapsto\mu(t,A)$:\[\{g_n\leq c\} = \bigcup_{A\in{\cal
P}_n}\mu(\cdot, A)^{-1}(-\infty, c\mathbb P(A)]\times A\in{\cal
T}\otimes{\cal S}_n\] Now clearly $\frac{d\mathbb
\mu(t,\cdot)}{d\mathbb P}|_{{\cal S}_n} = \mathbb E_{\mathbb
P}[\frac{d\mu(t,\cdot)}{d\mathbb P}|{\cal S}_n]$ a.s. , so that
$(g_n(t,\cdot))_n$ is a uniformly integrable $(\mathbb P, ({\cal
S}_n)_{n})$--martingale for each $t$. Hence there is $g(t,s)$ such
that
\[g_n(t,s)\to g(t,s)\qquad\text{$\mathbb P$--a.s. and in $L^1(\mathbb P)$, for each $t\in
T$}\]

\eproof

\end{document}